\documentclass[final]{article}

\usepackage{amsmath}
\usepackage{graphicx}
\usepackage{latexsym}
\usepackage{amssymb}
\usepackage{amsbsy}

\DeclareSymbolFont{msbm}{U}{msb}{m}{n}
\DeclareMathSymbol{\C}{\mathalpha}{msbm}{'103}
\DeclareMathSymbol{\R}{\mathalpha}{msbm}{'122}
\DeclareMathSymbol{\Z}{\mathalpha}{msbm}{'132}
\DeclareMathSymbol{\N}{\mathalpha}{msbm}{'116}

\newtheorem{remark}{Remark}
\newtheorem{proposition}{Proposition}
\newtheorem{theorem}{Theorem}
\newtheorem{definition}{Definition}

\def\RR{\mathbb R}

\def\be{\begin{equation}}
\def\ee{\end{equation}}
\def\bea{\begin{eqnarray}}
\def\ba{\begin{array}{l}\displaystyle}
\def\eea{\end{eqnarray}}
\def\ea{\end{array}}

\begin{document}
\title{Exponential Runge-Kutta methods for stiff kinetic equations}

\author{Giacomo Dimarco\thanks{Universit\'{e} de Toulouse; UPS, INSA, UT1, UTM;
CNRS, UMR 5219; Institut de Math\'{e}matiques de Toulouse; F-31062
Toulouse, France. ({\tt giacomo.dimarco@math.univ-toulouse.fr}).}
\and Lorenzo Pareschi\thanks{Mathematics Department, University of
Ferrara and CMCS, Ferrara, Italy ({\tt lorenzo.pareschi@unife.it}).}
}
\maketitle

\begin{abstract}
We introduce a class of exponential Runge-Kutta integration
methods for kinetic equations. The methods are based on a
decomposition of the collision operator into an equilibrium and a
non equilibrium part and are exact for relaxation operators of BGK
type. For Boltzmann type kinetic equations they work uniformly for
a wide range of relaxation times and avoid the solution of
nonlinear systems of equations even in stiff regimes. We give
sufficient conditions in order that such methods are
unconditionally asymptotically stable and asymptotic preserving.
Such stability properties are essential to guarantee the correct
asymptotic behavior for small relaxation times. The methods also
offer favorable properties such as nonnegativity of the solution
and entropy inequality. For this reason, as we will show, the
methods are suitable both for deterministic as well as
probabilistic numerical techniques.
\end{abstract}

\maketitle

{\bf Keywords:} Exponential integrators, Runge-Kutta methods, stiff equations,
Boltzmann equation, fluid limits, asymptotic preserving schemes.


\section{Introduction}
The numerical solution of the Boltzmann equation close to fluid
regimes represents a major computational challenge in rarefied gas
dynamics. In such regimes, in fact, the mean free path becomes very
small and standard computational approaches lose their efficiency
due to the necessity of using very small time steps in deterministic
schemes or, equivalently, a large number of collisions in
probabilistic approaches. Several authors have tackled the problem
in the past, and there is a large literature on the subject (see
\cite{BLM, CJR, degond1, toscani, Filbet, TK} and the references
therein). Most standard techniques are based on domain-decomposition
strategies and/or model reduction asymptotic methods. The direct
time discretization of the Boltzmann equation, in fact, represents a
challenge in such stiff regimes due to the high dimensionality and
the nonlinearity of the collision operator which makes unpractical
the use of implicit solvers. Exponential methods, like the so-called
time relaxed discretizations \cite{toscani}, combined with splitting
strategies represent one possible way to overcome the problem.
However, as discussed in various papers, the choice of the
Maxwellian equilibrium truncation in the schemes was based more on
physical considerations than on a direct mathematical derivation and
it is an open problem to determine an optimal truncation criteria
\cite{CaSa, PRHawaii}. Despite this, these discretizations have been
applied with success both in the context of spectral methods as well
as Monte Carlo methods \cite{Filbet2, PR}.

In this paper we propose a class of exponential integrators
\cite{EX1, EX2} for the homogeneous Boltzmann equation and related
kinetic equations which are based on explicit Runge-Kutta methods.
The main advantage of the approach here proposed is that it works
uniformly for a wide range of Knudsen numbers and avoids the
solution of nonlinear systems of equations even in stiff regimes.
Similarly to \cite{toscani}, we derive sufficient conditions such
that the resulting schemes can be represented as well as convex
combinations of density functions including a Maxwellian term.
This property is essential to achieve nonnegativity, physical
conservations and entropy inequality.

The starting point is the use of classical time spitting together
with a decomposition of the gain term of the collision operator
into an equilibrium and a non equilibrium part. Similar
decompositions for the distribution function has been used in
\cite{BLM} to derive unconditionally stable schemes and in
\cite{CPmc, dimarco1, dimarco2, dimarco3} to develop hybrid Monte
Carlo methods for kinetic equations. This decomposition of the
Boltzmann integral has been introduced recently by Jin and Filbet
\cite{Filbet} as a penalization method by a BGK-type relaxation
operator to derive Implicit-Explicit Runge-Kutta schemes that
overcome the stiffness of the full nonlinear collision operator.
We recall that analogous penalization techniques based on
linearized operators were previously used in the context of
exponential methods for parabolic equations with applications to
the Schr\"odinger equation \cite{EX1, EX2} and for kinetic
equations \cite{Lem, toscani}. In particular, in the present
paper, we will show that for space homogeneous equations the
Maxwellian truncation criteria introduced in \cite{toscani} is
equivalent to the penalization method in \cite{Filbet}. Let us
mention that a study of the accuracy property of splitting methods
in stiff regimes is beyond the scopes of the present paper.

Even if we develop our schemes using the Boltzmann equation as a
prototype model for its intrinsic difficulties and its relevance in
applications, the methods applies to any large system of stiff
ordinary differential equations of the form
\begin{equation}
Y'=\frac1{\varepsilon} R(Y),\quad Y(t_0)=Y_0,
\end{equation}
where $\varepsilon>0$ is a small parameter, $Y\in\R^N$ and the
non-linear operator $R(Y)$ is a dissipative relaxation operator as
in \cite{CLL}. Such operator is endowed with a $n\times N$ matrix
$Q$ of rank $n<N$ such that $QR(Y)=0$, $\forall\,\, Y$. This gives
a vector of $n$ conserved quantities $y=QY$. Solutions which
belong to the kernel of the operator $R(Y)=0$ are uniquely
determined by the conserved quantities $Y=E(y)$ and characterize
the manifold of local equilibria. Important examples of such
dissipative relaxation operators arise in discrete kinetic theory,
shallow water equations, granular gases, traffic flows and in
general in finite difference/volume discretizations of several
kinetic equations.

The methods here proposed use the following decomposition
\cite{Filbet}
\begin{equation}
R(Y)=N(Y)+A(E(y)-Y),
\end{equation}
where $N(Y)$ represents the non-dissipative non-linear part, $A>0$
is an estimate of the Jacobian of $R$ evaluated at equilibrium and
$E(y)-Y$ is a simple dissipative linear relaxation operator. Note
that, at variance with standard linearization techniques which
operate on the short time scale, the operator is linearized on the
asymptotically large time scale.

This decomposition permits to apply exponential techniques which
solve exactly the linear part and are explicit in the non-linear
part. The use of such techniques, as we will see, is essential in
order to achieve some unconditional stability properties of the
numerical schemes \cite{MZ}. Such properties are usually
characterized as unconditional asymptotic stability and asymptotic
preservation.

Here we derive sufficient conditions for asymptotic stability and
asymptotic preservation. This permits to introduce asymptotically
stable and asymptotic preserving methods up to order 4 which are
exact for BGK-type kinetic equations.

The rest of the paper is organized as follows. First in Section 2 we
present the Boltzmann equation and its fluid-dynamic limit. Operator
splitting methods and various notions of stability are illustrated
in Section 3. Next, in Section 4, we introduce the explicit
exponential Runge-Kutta methods and derive conditions for asymptotic
stability and asymptotic preservation. Examples of methods up to
order 4 are also constructed. In Section 5 we describe numerical
experiments for homogeneous problems and one application to the full
Boltzmann equation. In a final section we discuss conclusions and
perspectives of the methods proposed in this article.

\section{The Boltzmann equation}

We consider the Boltzmann equation of rarefied gas dynamics \cite{cercignani}
\be
\partial_t f + v\cdot\nabla_{x}f
=Q(f,f)\label{eq:1}\ee with initial data \be
f|_{t=0}=f_{0}.\label{eq:ini}\ee Here $f(x,v,t)$ is a non negative
function describing the time evolution of the distribution of
particles with velocity $v \in \R^{3}$ and position $x \in \Omega
\subset \R^{d_x}$ at time $ t
> 0$. In the sequel for notation simplicity we will omit the
dependence of $f$ from the independent variables $x,v,t$ unless
strictly necessary. The operator $Q(f,f)$ which describes
particles interactions has the form \be Q(f,f)=\int_{\RR^3\times
S^2} B(|v-v_*|,n) [f(v')f(v'_*)-f(v)f(v_*)]\,dv_*\,dn \label{eq:Q}
\ee where \be v'=v+\frac12(v-v_*)+\frac12|v-v_*|n,\quad
v'_*=v+\frac12(v-v_*)-\frac12|v-v_*|n, \ee and $B(|v-v_*|,n)$ is a
nonnegative collision kernel characterizing the microscopic
details of the collision given by \[
B(|v-v_*|,n)=\sigma\left(\frac{(v-v_*)}{|v-v_*|}\cdot n\right)
|v-v_*|^\gamma,
\] with $\gamma \in [0,3)$. The case $\gamma=1$ is
referred to as hard spheres case, whereas the simplified situation
$\gamma=0$, is referred to as Maxwell case. Note that in most
applications the angle dependence is ignored and $\sigma$ is
assumed constant \cite{bird}.

The operator $Q(f,f)$ is such that the the local conservation
properties are satisfied \be\int_{\R^{3}} mQ(f,f) dv=:\langle
mQ(f,f)\rangle=0 \label{eq:QC}\ee where $m(v)=(1,v,\frac{|v|^2}{2})$
are the collision invariants. In addition it satisfies the entropy
inequality \be \frac{d}{dt}\int_{\R^{3}}f\log f\,dv = \int_{\R^{3}}
Q(f,f)\log f dv \leq 0. \label{eq:entropy} \ee Integrating
(\ref{eq:1}) against its invariants in the velocity space leads to
the following set of non closed conservations laws \be
\partial_t \langle mf\rangle+\nabla_x
\langle vmf\rangle=0.\label{eq:macr}\ee Equilibrium functions for
the operator $Q(f,f)$ (i.e. solutions of $Q(f,f)=0$) are local
Maxwellian of the form \be M_{f}(\rho,u,T)=\frac{\rho}{(2\pi
T)^{3/2}}\exp\left(\frac{-|u-v|^{2}}{2T}\right), \label{eq:M}\ee
where $\rho$, $u$, $T$ are the density, mean velocity and
temperature of the gas in the x-position and at time $t$ defined
as \be (\rho,\rho u,E)=\langle m f \rangle, \qquad
T=\frac1{3\rho}(E-\rho|u|^2). \ee We will denote by \be
U=(\rho,u,E),\qquad M[U]=M_f, \ee clearly we have \be U=\langle
mM[U]\rangle.\ee Now, when the mean free path between particles is
very small the operator $Q(f,f)$ is large and we can rescale the
space and time variables in (\ref{eq:1}) as \be x'=\varepsilon x,
\ \ t'=\varepsilon t\ee to obtain \be
\partial_t f + v\cdot\nabla_{x}f
=\frac1{\varepsilon}Q(f,f)\label{eq:1b}\ee where $\varepsilon$ is
a small parameter proportional to the mean free path and the
primes have been omitted to keep notation simple.

Passing to the limit for $\varepsilon\rightarrow 0$ leads to
$f\rightarrow M[U]$ and thus we obtain the closed hyperbolic system
of compressible Euler equations for the macroscopic variables $U$
\be
\partial_t U+\nabla_xF(U)=0
\label{eq:Euler} \ee with
\[
F(U)=\langle vmM[U]\rangle=(\rho u, \varrho u \otimes u+pI,
Eu+pu),\quad p=\rho T,
\]
where $I$ is the identity matrix.

For small but non zero values of the Knudsen number, the evolution
equation for the moments can be derived by the so-called
Chapman-Enskog expansion \cite{cercignani}. This approach gives
the compressible Navier-Stokes equations as a second order
approximation with respect to $\varepsilon$ to the solution to the
Boltzmann equation.

\subsection{Operator splitting and stability definitions}

Here we restrict ourselves to operator splitting based schemes. It
is well-known, in fact, that most numerical methods for the
Boltzmann equation are based on a splitting in time between free
particle transport and collisions \cite{bird, Filbet2}. Possible
extension of the present theory to non-splitting schemes is actually
under study and it will be considered elsewhere.

Even if it is difficult to give a rigorous definition of asymptotic
preserving scheme since the concept has been used for a long time
in the physics and mathematics literature and may refer to different
discretization parameters, here
following \cite{Jin, Jin2, PRimex} we formalize this notion for the time
discretization of equation
(\ref{eq:1}).
\begin{definition}
A consistent time discretization method for (\ref{eq:1}) of
stepsize $\Delta t$ is {\em asymptotic preserving (AP)} if,
independently of the initial data (\ref{eq:ini}) and of the
stepsize $\Delta t$, in the limit $\varepsilon\to 0$ becomes a
consistent time discretization method for the reduced system
(\ref{eq:Euler}).
\end{definition}

Note that this definition does not imply that the scheme preserves
the order of accuracy in $t$ in the stiff limit $\varepsilon\to
0$. In the latter case the scheme is said {\em asymptotically
accurate}.

As discussed above, the starting point in the solution of the
kinetic equations is given by an operator splitting of
(\ref{eq:1}) in a time interval $[0,\Delta t]$ between relaxation
\be
\partial_t f=\frac{1}{\varepsilon}Q(f,f),\label{eq:9}\ee
and free transport \be
\partial_t f + v\cdot\nabla_{x}f
=0. \label{eq:tr} \ee This situation is typical of Monte Carlo
methods and of several other numerical codes used in realistic
simulations. Even if this splitting, usually referred to as
Lie-Trotter splitting, is limited to first order it permits to
treat separately the hyperbolic free transport from the stiff
relaxation step which often is of paramount importance in
applications.

Higher order splitting formulas can be derived in different ways
(see \cite{HLW}). Let us denote by ${\cal T}_{\Delta t}(f)$ and
${\cal C}_{\Delta t}(f)$ the above transport and collision steps
in a time interval $[0,\Delta t]$ starting from the initial data
$f_0$, then the well-known second order Strang splitting
\cite{Strang} can be written as \be {\cal C}_{\Delta t/2}({\cal
T}_{\Delta t}({\cal C}_{\Delta t/2}(f_0))). \label{eq:s2} \ee

Unfortunately for splitting methods of order higher then two it can
be shown that it is impossible to avoid negative time steps both in
the transport as well as in the collision \cite{HLW}. Higher order
formulas which avoid negative time stepping can be obtained as
suitable combination of splitting steps \cite{DS}. Note, however,
that the appearance of negative coefficients or negative time steps
in high order formulas may lead to several drawbacks in practical
applications like the lack of positivity of the solution which makes
very difficult their use in Monte Carlo schemes.

As mentioned before the study of the accuracy properties of the
different splitting methods in stiff regimes, although important,
is beyond the scopes of the present manuscript and we refer to
\cite{JL, Th} for an error analysis of some splitting methods in
stiff conditions.

Now we can reformulate the asymptotic preserving property and
prove that

\begin{proposition}
A sufficient condition for a consistent time discretization method
of stepsize $\Delta t$ applied to the operator splitting
approximation of (\ref{eq:1}), given by
(\ref{eq:9})-(\ref{eq:tr}), to be AP is that the time
discretization of step (\ref{eq:9}), independently of the initial
data (\ref{eq:ini}) and of the stepsize $\Delta t$, in the limit
$\varepsilon\to 0$ projects the solution $f$ over the local
Maxwellian equilibrium $M[U_0]$, $U_0=\langle mf_0\rangle$.
\label{pr:1}
\end{proposition}
The proof of the above proposition is an immediate consequence of
the fact that as $\varepsilon\to 0$ step (\ref{eq:9}) degenerates
into the projection ${\cal C}_{\Delta t}(f_0)=M[U_0]$ which
coupled with the transport step (\ref{eq:tr}) originates a
so-called kinetic approximation \cite{CP} to the Euler equation
(\ref{eq:Euler}) given by ${\cal T}_{\Delta t}(M[U_0])$. We omit
further details.

In other words, Proposition \ref{pr:1} states that if the
relaxation step (\ref{eq:9}) is AP then the whole splitting
(\ref{eq:9})-(\ref{eq:tr}) is AP. Analogous results hold true for
higher order splitting methods.

In the sequel we will focus on the solution to the space
homogeneous Boltzmann equation (\ref{eq:9}). In fact, most
computational challenges related to the behavior of the full
equation for small values of $\varepsilon$ depend on the time
discretization of the homogeneous step.

Of course AP is an important property in term of stability of the
numerical scheme in stiff regimes. For implicit Runge-Kutta
methods applied to (\ref{eq:9}) it has been shown in \cite{PRimex}
that AP is equivalent to the notion of $L$-stability
\cite{HairerWanner}.

For general unconditionally stable schemes a weaker requirement is
the notion of {\em asymptotic stability} (AS) \cite{MZ}. Let us
denote by $f(t)$ and $g(t)$ two solutions of (\ref{eq:9})
corresponding respectively to the initial data $f_0$ and $g_0$
such that $U_0=\langle m f_0 \rangle=\langle m g_0 \rangle$. It
can be proved that, for a suitable distance $d(\cdot,\cdot)$,
system (\ref{eq:9}) is contractive $d(f(t),g(t)) \leq d(f_0,g_0)$,
and asymptotically stable since $f(t)\to M[U_0]$ and $g(t)\to
M[U_0]$ as $t\to\infty$ (see \cite{cercignani} for example) and so
$f(t)-g(t) \to 0$ as $t\to\infty$. We refer to \cite{TV} for more
details and examples of contractive metrics for problem
(\ref{eq:9}) in the case of Maxwell molecules. A particular metric
is presented in Section \ref{sec:3-3}.

Let us denote by $f^n$ and $g^n$, $n\geq 1$ the numerical solution
at $t=n\Delta t$ obtained with a given time discretization method
applied to (\ref{eq:9}) with initial data $f_0$ and $g_0$
respectively. Now we can introduce the following definition.

\begin{definition} A time discretization method for (\ref{eq:9})
is called {\em unconditionally contractive} with respect to the
distance $d(\cdot,\cdot)$ if $d(f^1,g^1) \leq d(f_0,g_0)$ holds
for all $f_0$, $g_0$ such that $\langle m f_0 \rangle=\langle m
g_0 \rangle$ and for all stepsizes $\Delta t$. Furthermore, it is
called {\em unconditionally asymptotically stable} if $f^n-g^n\to
0$ as $n\to\infty$ independently of the step size $\Delta t$.
\end{definition}

Note that unlike contractivity, asymptotic stability is not
related to a specific metric. Contractivity of Runge-Kutta methods
has been studied in \cite{Kr} and it is well-known that such
methods have limited order of accuracy. Clearly for implicit
Runge-Kutta methods, asymptotic stability is closely related to
$A$-stability.



For the sake of completeness we finally introduce the notion of
entropic stability, namely schemes that preserve the entropy
inequality (\ref{eq:entropy}).

\begin{definition}
A time discretization method for (\ref{eq:9}) is called {\em
unconditionally entropic} if $H(f^{n+1})\leq H(f^n)$, where
$H(f)=\int_{\R^{3}} f\log f\,dv$, independently of the step size
$\Delta t$.
\end{definition}

The above monotonicity property for Runge-Kutta schemes is
essentially equivalent to the so-called {strong stability
preserving} (SSP) property which is often used when dealing with
the time discretization of partial differential equations
\cite{GST}. Let us recall that Implicit Euler is the sole
unconditional SSP method. All SSP Runge-Kutta methods of order
greater than one suffer from some time-step restriction \cite{Hi}.
Thus, except for first order implicit Euler, the entropy
inequality is not satisfied by high order unconditionally stable
Runge-Kutta schemes applied to (\ref{eq:9}) unless a suitable time
step restriction is considered.

\section{Exponential methods}

Since we aim at developing unconditionally stable schemes, the most
natural choice would be to use implicit solvers applied to
(\ref{eq:9}). Unfortunately the use of fully implicit schemes for
(\ref{eq:9}) is unpracticable due to the prohibitive computational
cost required by the solution of the very large non-linear algebraic
system originated by the five fold integral appearing in $Q(f,f)$
which has to be computed in each spatial cell at each time step in
the inhomogeneous cases. We will see in this section a possible way
to overcome these difficulties.

First we rewrite the homogeneous equation (\ref{eq:9}) in the form
\be
\partial_t f=\frac{1}{\varepsilon}(P(f,f)-\mu f),\label{eq:10}
\ee where $P(f,f)=Q(f,f)+\mu f$ and $\mu>0$ is a constant such
that $P(f,f)\geq 0$. Typically $\mu$ is an estimate of the largest
spectrum of the loss part of the collision operator. Let us
emphasize that most of the subsequent theory can be generalized to
the case where $P(f,f)$ is not strictly positive and $\mu$ is an
arbitrary nonnegative constant. However such an assumption makes
the presentation easier and we will discuss possible
generalizations later.

By construction we have the following \be \frac1{\mu}\langle m
P(f,f)\rangle = \langle m f\rangle=U. \ee Thus $P(f,f)/\mu$ is a
density function and we can consider the following
decomposition \be P(f,f)/\mu=M[U]+g.\label{eq:dec}\ee The function $g$
represents the non equilibrium part of the distribution function
and from the definition above it follows that $g$ is in general
non positive. Moreover since $P(f,f)/\mu$ and $M[U]$ have the same
moments we have \be \langle mg\rangle=0.\ee

The homogeneous equation can be written in the form
\be
\partial_t
f=\frac{\mu}{\varepsilon}g+\frac{\mu}{\varepsilon}(M[U]-f)=\frac{\mu}{\varepsilon}\left(\frac{P(f,f)}{\mu}-M[U]\right)+
\frac{\mu}{\varepsilon}(M[U]-f).\label{eq:11}\ee The above system
is equivalent to the penalization method for the collision
operator recently introduced in \cite{Filbet}. Note that even if
$M[U]$ is nonlinear in $f$, thanks to the conservation properties
(\ref{eq:QC}), it remains constant during the relaxation process.
The main feature of such formulation is that on the right hand
side we have a stiff dissipative linear part $\mu
(M[U]-f)/\varepsilon$ which characterizes the asymptotic behavior
of $f$ and a stiff non dissipative non linear part $(P(f,f)/\mu-
M[U])/\varepsilon$ which describes the deviations of $P(f,f)/\mu$
from $M[U]$, or equivalently the deviations of the Boltzmann
operator from a BGK-like relaxation term.

We remark that the decision whether problem (\ref{eq:9}), or
equivalently (\ref{eq:10}) and (\ref{eq:11}), should be regarded
as stiff or nonstiff does not depend only on the ratio
$\mu/\varepsilon$ but depends also on the chosen initial
conditions. If the initial data is close to local equilibrium
$f=M[U]+O(\varepsilon)$, then the problem is clearly nonstiff. In
fact, if $f=M[U]+\varepsilon f_1$ we have
\[
Q(f,f)=\varepsilon [Q(f_1,f)+Q(f,f_1)+\varepsilon Q(f_1,f_1)],
\]
and so $Q(f,f)=O(\varepsilon)$ and there is no need of using a
specific stiff solver. For this reason, and the fact that the
transport step (\ref{eq:tr}) drives the solution far from local
equilibrium, in the sequel we concentrate our analysis to non
equilibrium initial data.

In such case the problem is stiff as a whole and a fully implicit
method should be used in the numerical integration to avoid
stability constraints of the type $\Delta t = O(\varepsilon)$. On
the other hand the linear part itself suffices to characterize the
correct large time behavior of $f$. Therefore, instead of fully
implicit methods, one should hopefully use methods which are
implicit in the linear part and explicit in the non-linear part. One
class of such methods is given by the IMEX Runge-Kutta schemes
\cite{Ascher, PRimex}. Note, however, that standard IMEX schemes may
lose their stability properties since here also the explicit part is
stiff. An alternative approach is based on the so-called exponential
integrators where the exact solution of the linear part is used in
the construction of the numerical methods \cite{EX1, EX2}.

\subsection{Explicit exponential Runge-Kutta schemes}

In order to derive the methods it is useful to rewrite
(\ref{eq:11}) as \be \frac{\partial (f-M)e^{\mu
{t}/{\varepsilon}}}{\partial t}=\frac1{\varepsilon}(P(f,f)-\mu
M)e^{\mu {t}/{\varepsilon}}. \label{eq:15} \ee

The above form is readily obtained if one multiplies (\ref{eq:11})
by the integrating factor $\exp(\mu t/\varepsilon)$ and takes into
account the fact that $M$ does not depend of time. A class of
explicit exponential Runge-Kutta schemes is then obtained by
direct application of an explicit Runge-Kutta method to
(\ref{eq:15}). More in general we can consider the family of
methods characterized by
\begin{eqnarray}
\nonumber F^{(i)}&=&e^{-c_i\mu\Delta
t/\varepsilon}f^n+\frac{\mu\Delta t}{\varepsilon}
\sum_{j=1}^{i-1}A_{ij}(\mu\Delta
t/\varepsilon)\left(\frac{P(F^{(j)},F^{(j)})}{\mu}-M^n\right)\\
\label{eq:rk1}
\\[-.25cm]
\nonumber &+& \left(1-e^{-c_i\mu\Delta t/\varepsilon}\right)M^n,
\qquad
i=1,\ldots,\nu\\
\nonumber f^{n+1}&=&e^{-\mu\Delta
t/\varepsilon}f^n+\frac{\mu\Delta
t}{\varepsilon}\sum_{i=1}^{\nu}W_i(\mu\Delta
t/\varepsilon)\left(\frac{P(F^{(i)},F^{(i)})}{\mu}-M^n\right)\\
\label{eq:rk2}
\\[-.25cm]
\nonumber &+&\left(1-e^{-\mu\Delta t/\varepsilon}\right)M^n,
\end{eqnarray}
where $\Delta t$ is the time step, $f^n=f(t^n)$, $M^n=M(t^n)$,
$c_i \geq 0$, and the coefficients $A_{ij}$ and the weights $W_i$
are such that
\[
A_{ij}(0)=a_{ij},\quad W_i(0)=w_i,\quad i,j=1,\ldots,\nu
\]
with coefficients $a_{ij}$ and weights $w_i$ given by a standard
explicit Runge-Kutta method called the underlying method. Various
schemes come from the different choices of the underlying method.
The most popular approaches are the integrating factor (IF) and
the exponential time differencing (ETD) methods \cite{EX1, EX2,
MZ}. Since $M^n$ does not depend on time during the collision
process in the sequel we will omit the index $n$.

For the so-called Integrating Factor methods, which correspond to
a direct application of the underlying method to (\ref{eq:15}), we
have
\begin{eqnarray}
\nonumber
A_{ij}(\lambda)&=&a_{ij}e^{-(c_i-c_j)\lambda},\quad
i,j=1,\ldots,\nu,\quad i > j\\[-.1cm]
\\
\nonumber
W_i(\lambda)&=&w_i e^{-(1-c_i)\lambda},\quad i=1,\ldots,\nu,
\end{eqnarray}
with $\lambda=\mu\Delta t/\varepsilon$.

The first order IF scheme reads
\begin{equation}
f^{n+1}=e^{-\frac{\mu \Delta t}{\varepsilon}}
f^{n}+\frac{\mu\Delta t}{\varepsilon}e^{-\frac{\mu \Delta
t}{\varepsilon}}\left(\frac{P(f^n,f^n)}{\mu}-M\right)+\left(1-e^{-\frac{\mu
\Delta t}{\varepsilon}}\right)M, \label{eq:if1b}
\end{equation}
which is based on explicit Euler. For such methods the order of
accuracy is the same as the order of the underlying method.

The Exponential Time Differencing methods are strictly connected
with the integral representation of (\ref{eq:15}). In the general
case the coefficients for ETD methods have the form
\begin{eqnarray*}
A_{ij} (\lambda) &=& \int_0^1 e^{(1-s)c_i\lambda} p_{ij}(s)\,ds,
\quad
i,j=1,\ldots,\nu,\quad i > j\\
W_{i} (\lambda) &=& \int_0^1 e^{(1-s)\lambda} p_{i}(s)\,ds, \quad
i=1,\ldots,\nu,
\end{eqnarray*}
where $p_i$ and $p_{ij}$ are suitable polynomials.

The standard first order ETD method based on explicit Euler in our
case gives \cite{EX1}
\begin{equation}
f^{n+1}=e^{-\frac{\mu \Delta t}{\varepsilon}}
f^{n}+\frac{\mu\Delta t}{\varepsilon}\varphi\left(\frac{\mu\Delta
t}{\varepsilon}\right)\frac{P(f^n,f^n)}{\mu}, \label{eq:etd1c}
\end{equation}
where $\varphi(z)=(1-e^{-z})/{z}$.

\subsection{Time Relaxed methods}
A class of exponential methods for kinetic equations, the
so-called time relaxed (TR) methods, has been introduced in
\cite{toscani} as a combination of an exponential expansion (or
Wild sum) together with a suitable Maxwellian truncation. In this
paragraph we show that these schemes included already
decomposition (\ref{eq:11}) and can be derived directly from a
suitable Taylor expansion of (\ref{eq:15}).

To show this, let us first introduce the change of variables
\[
\tau=1-\exp(-\mu t/\varepsilon),
\]
which, using the bilinearity of $P(f,f)$, gives the equation \be
\frac{\partial}{\partial \tau}\left[
(f-M)\frac{1}{1-\tau}\right]=(P(f,f)-\mu
M)\frac{1}{\mu(1-\tau)^2}. \label{eq:16} \ee The application of an
explicit Runge-Kutta scheme to (\ref{eq:16}) with time step
$\Delta \tau=1-\exp(-\mu\Delta t/\varepsilon)$ leads to a class of ETD methods.
For example the first order scheme based
on explicit Euler in the original variables yields
\begin{equation}
f^{n+1}=e^{-\frac{\mu \Delta t}{\varepsilon}}
f^{n}+\frac{\mu\Delta
t}{\varepsilon}\varphi_1\left(\frac{\mu\Delta
t}{\varepsilon}\right)\left(\frac{P(f^n,f^n)}{\mu}-M\right)+(1-e^{-\frac{\mu
\Delta t}{\varepsilon}}) M, \label{eq:etd1b}
\end{equation}
where $\varphi_k(z)=e^{-z}({1-e^{-z}})^k/{z}$, $k=1,2,\ldots$.

Note that such scheme coincides with the first order exponential
time relaxed method derived in \cite{toscani} and differs from the
standard ETD method based on explicit Euler. Higher order ETD
methods can be derived as well simply applying higher order
explicit Runge-Kutta methods to (\ref{eq:16}). Although
interesting, here we do not explore further this class of schemes.

Now let us consider a different approach by taking the Taylor
expansion of $(f-M)/(1-\tau)$ around $\tau=0$ in (\ref{eq:16}).
This leads to
\begin{eqnarray*}
(f-M)/(1-\tau)&=&(f_0-M)+\tau\left[\frac{P(f_0,f_0)}{\mu}-M\right]\\
&+&\frac{\tau^2}{2}\left[\frac{P(P(f_0,f_0),f_0)+P(f_0,P(f_0,f_0))}{\mu^2}-2M\right]+O(\tau^3)
\end{eqnarray*}
where we have used the bilinearity of the operator $P(f,f)$.

If we compute all the terms in the expansion and use recursively
the bilinearity of $P(f,f)$ we can state the following

\begin{proposition}
The solution to problem (\ref{eq:11}) or equivalently
(\ref{eq:15}) or (\ref{eq:16}) can be represented as
 \be f(v,t) = (1-\tau) f_0(v)+
(1-\tau)\sum_{k=1}^{\infty} \tau^k (f_k^n(v)-M(v)) + \tau M(v),
\label{eq:ws2}\ee where the functions $f_{k}$ are given by the
recurrence formula \be f_{k+1}(v) = \frac{1}{k+1}\sum_{h=0}^k
\frac{1}{\mu} P ( f_{h}, f_{k-h})(v),\quad k=0,1,\ldots.
\label{eq:CF-r} \ee
\end{proposition}

By truncating expansion (\ref{eq:ws2}) at the order $m$, and
reverting to the old variables, we recover exactly the time
relaxed schemes presented in \cite{toscani} \be f^{n+1} =
e^{-\mu{\Delta t}/{\varepsilon}}f^n+e^{-\mu{\Delta
t}/{\varepsilon}}\sum_{k=1}^m (1-e^{-\mu{\Delta
t}/{\varepsilon}})^k (f_k^n-M) + (1-e^{-\mu{\Delta
t}/{\varepsilon}}) M, \ee which, using the fact that
\[
1-e^{-\mu{\Delta t}/{\varepsilon}}\sum_{k=0}^m (1-e^{-\mu{\Delta
t}/{\varepsilon}})^k = (1-e^{-\mu{\Delta
t}/{\varepsilon}})^{m+1},
\]
can be rewritten in the usual form
emphasizing their convexity properties
 \be f^{n+1} = e^{-\mu{\Delta
t}/{\varepsilon}}\sum_{k=0}^m (1-e^{-\mu{\Delta
t}/{\varepsilon}})^k f_k^n + (1-e^{-\mu{\Delta
t}/{\varepsilon}})^{m+1} M. \label{eq:trs2}\ee A remarkable
feature of these methods is that the functions $f_k(v)$ are
density functions with the same moments of the initial data $
\langle m f_k \rangle=\langle f_0 \rangle. $ Such property,
together with unconditional nonnegativity and convexity of the
weights, is enough to guarantee asymptotic preservation of the
schemes as well as nonnegativity, contractivity and entropic
stability (see \cite{toscani} for details).

Clearly TR schemes do not belong to the general family of methods
defined by (\ref{eq:rk1})-(\ref{eq:rk2}), since they are based on
the assumption that $P(f,f)$ is a bilinear operator.

\subsection{Main properties}
\label{sec:3-3} In this section we will establish the main
properties for a general exponential scheme in the form
(\ref{eq:rk1})-(\ref{eq:rk2}). For some of the properties, like
contractivity and entropic stability, we give proofs in the case of
Maxwellian molecules.

Since solutions to (\ref{eq:9}) are nonnegative and preserves the
mass in our analysis we will restrict, without loss of generality,
to probability density functions.

Let us denote by ${\cal P}_2(\R^{3})$, the class of all
probability density functions $f$ on $\R^{3}$, such that \be
\int_{\R^{3}} |v|^2 dF(v) < \infty. \ee We introduce a metric on
${\cal P}_2(\R^{3})$ by \be d_2(f,g) = \sup_{\xi\in
\R^{3}}\frac{|\hat f(\xi)-\hat g(\xi)|}{|\xi|^2} \ee where $\hat
f$ is the Fourier transform of $f$ \be \hat f(\xi) = \int_{\R^d}
e^{-i\xi\cdot v} dF(v). \ee

The metric $d_2(\cdot,\cdot)$ is nonexpanding with time along
trajectories of the Boltzmann equation; that is, if $f$ and $g$ are
two such solutions $d_2(f(t), g(t))\leq d_2(f_0, g_0)$. The above
property is a consequence of the fact that for Maxwell's molecules
we have \be Q(f,f)=\int_{\RR^3\times S^2}
\sigma\left(\frac{v-v_*}{|v-v_*|}\cdot n\right)
[f(v')f(v'_*)-f(v)f(v_*)]\,dv_*\,dn \label{eq:QM} \ee where for any
fixed unit vector ${\bar e}$
\[
\int_{S^2}\sigma({\bar e}\cdot n)\,dn=S,
\]
with $S>0$ a constant. Taking $\mu=S$ and
\[
P(f,f)=\int_{\RR^3\times S^2}
\sigma\left(\frac{v-v_*}{|v-v_*|}\cdot n\right)
f(v')f(v'_*)\,dv_*\,dn,
\]
we have $Q(f,f)=P(f,f)+\mu f$ and \cite{TV} \be
d_2(P(f,f),P(g,g))\leq \mu d_2(f,g). \ee We refer to \cite{TV} for
more information on the above metric and other contractive metrics
for the Boltzmann equation in the case of Maxwellian molecules.

Now let us denote by $f^{n}$ and $g^n$ the corresponding solutions
obtained with an explicit exponential Runge-Kutta method. We have
\begin{eqnarray}
\nonumber F^{(i)}-G^{(i)}&=&e^{-c_i\mu\Delta
t/\varepsilon}(f^n-g^n)+\frac{\Delta t}{\varepsilon}
\sum_{j=1}^{i-1}A_{ij}(\mu\Delta
t/\varepsilon)(P(F^{(j)},F^{(j)})-P(G^{(j)},G^{(j)}))\\
\nonumber f^{n+1}-g^{n+1}&=&e^{-\mu\Delta
t/\varepsilon}(f^n-g^n)+\frac{\Delta
t}{\varepsilon}\sum_{i=1}^{\nu}W_i(\mu\Delta
t/\varepsilon)(P(F^{(i)},F^{(i)})-P(G^{(i)},G^{(i)})),
\end{eqnarray}
and then
\begin{eqnarray}
\nonumber d_2(F^{(i)},G^{(i)})&\leq &e^{-c_i\mu\Delta
t/\varepsilon}d_2(f^n,g^n)+\frac{\mu \Delta t}{\varepsilon}
\sum_{j=1}^{i-1}|A_{ij}(\mu\Delta
t/\varepsilon)|d_2(F^{(j)},G^{(j)})\\
\nonumber d_2(f^{n+1},g^{n+1})&\leq &e^{-\mu\Delta
t/\varepsilon}d_2(f^n,g^n)+\frac{\mu \Delta
t}{\varepsilon}\sum_{i=1}^{\nu}|W_i(\mu\Delta
t/\varepsilon)|d_2(F^{(i)},G^{(i)}).
\end{eqnarray}
Let us denote by $\lambda=\mu\Delta t/\varepsilon$, $\bar
A(\lambda)$ the $\nu\times\nu$ matrix of elements
$|A_{ij}(\lambda)|$, $\bar w(\lambda)$ the $\nu\times 1$ vector of
elements $|W_i(\lambda)|$, $\bar d$ the $\nu\times 1$ vector of
elements $d_2(F^{(i)},G^{(i)})$ and $\bar e$ the $\nu\times 1$
unit vector we can write
\begin{eqnarray*}
\nonumber \left(I-\lambda {\bar A}(\lambda)\right){\bar d}&\leq
&{\bar E}(\lambda) d_2(f^n,g^n){\bar e}\\
\nonumber d_2(f^{n+1},g^{n+1})&\leq
&e^{-\lambda}d_2(f^n,g^n)+\lambda{\bar w}(\lambda)^T{\bar d},
\end{eqnarray*}
where ${\bar
E}(\lambda)={\rm diag}(e^{-c_1\lambda},\ldots,e^{-c_\nu\lambda})$.

Then we obtain
\be d_2(f^{n+1},g^{n+1})\leq
R(\lambda)d_2(f^{n},g^{n}),\ee where \begin{eqnarray}
R(\lambda)&=&e^{-\lambda}+\lambda {\bar
w}(\lambda)^T(I-\lambda{\bar A}(\lambda))^{-1}{\bar E}(\lambda){\bar e}\\
&=& e^{-\lambda}+\sum_{k=0}^{\nu-1}\lambda^{k+1} {\bar
w}(\lambda)^T{\bar A}(\lambda)^k{\bar E}(\lambda){\bar e}.
\label{eq:R}\end{eqnarray} To derive the last expression we expanded
$(I-\lambda{\bar A}(\lambda))^{-1}$ by the geometric series and use
the fact that ${\bar A}(\lambda)$ is strictly lower triangular and
so it is a nilpotent matrix of degree $\nu$.

We can state \cite{MZ}
\begin{theorem}
If an explicit exponential Runge-Kutta method in the form
(\ref{eq:rk1})-(\ref{eq:rk2}) satisfies \be R(\lambda)\leq 1,\quad
\forall\,\,\lambda \geq 0, \label{eq:cond}\ee with $R(\lambda)$
given by (\ref{eq:R}) then it is unconditionally contractive and
unconditionally stable with respect to the metric
$d_2(\cdot,\cdot)$.
\end{theorem}

Note that for an IF method we have
\[
|A_{ij}(\lambda)|\leq |a_{ij}|e^{-(c_i-c_j)\lambda},\quad
|W_i(\lambda)|\leq |w_{i}|e^{-(1-c_i)\lambda},
\]
thus we require \be 0=c_1\leq c_2\ldots\leq c_\nu\leq 1,
\label{eq:cond2}\ee in order for the above quantities to be
bounded independently of $\lambda$.

We have \be
R(\lambda)=e^{-\lambda}\left(1+\sum_{k=0}^{\nu-1}\lambda^{k+1}
{\bar w}^T{\bar A}^k{\bar e}\right). \label{eq:RIF}\ee Thus
condition (\ref{eq:cond}) is satisfied when \be
1+\sum_{k=0}^{\nu-1}\lambda^{k+1} {\bar w}^T{\bar A}^k{\bar e}\leq
e^\lambda=1+\sum_{k=0}^\infty \frac{\lambda^{k+1}}{(k+1)!}.
\label{eq:n} \ee Now, if we consider an underlying Runge-Kutta
method with a $\nu\times\nu$ non negative coefficient matrix $A$
and a $\nu\times 1$ non negative vector of weights $w$ then
(\ref{eq:n}) holds if
\[
{w}^T{A}^k{\bar e}\leq \frac{1}{(k+1)!},\quad k=0,1,\ldots,\nu-1.
\]
The above condition is clearly satisfied if the underlying
Runge-Kutta method is a $\nu$-stage explicit Runge-Kutta method of
order $\nu$.

Thus we have proved
\begin{proposition}
An explicit IF method is unconditionally contractive and
asymptotically stable with respect to the metric
$d_2(\cdot,\cdot)$ if the underlying Runge-Kutta method is a
$\nu$-stage explicit Runge-Kutta method of order $\nu$ with
nonnegative coefficients and weights satisfying (\ref{eq:cond2}).
\end{proposition}

As pointed out in \cite{Kr, MZ} examples of such methods are
well-known up to $\nu = 4$ and the classical RK method of order four
is the sole method with four stages. Moreover, it is also known that
there does not exist explicit methods of order greater than four
satisfying (\ref{eq:cond2}) with positive weights. For methods which
are not of IF type there are no accuracy barrier, for example all
time relaxed method satisfy immediately condition (\ref{eq:cond}).
Other examples are reported in \cite{EX2, MZ}.

We have

\begin{theorem}
If an explicit exponential Runge-Kutta method in the form
(\ref{eq:rk1})-(\ref{eq:rk2}) satisfies \be \lim_{\lambda \to
\infty}R(\lambda)= 0, \label{eq:condAP}\ee with $R(\lambda)$ given
by (\ref{eq:R}) then it is asymptotic preserving.
\end{theorem}

In fact taking $g_0=M$ we have
\[
d_2(f^{n+1},M)\leq R(\lambda) d_2(f^n,M),
\]
and so $d_2(f^{n+1},M)$ goes to $0$ as $\lambda\to\infty$.

For an IF method $R(\lambda)$ is given by (\ref{eq:RIF}) and
condition (\ref{eq:condAP}) is always satisfied. Thus

\begin{proposition}
An explicit IF method is asymptotic preserving if the underlying
Runge-Kutta method satisfies (\ref{eq:cond2}). \label{Pr:4}
\end{proposition}

For practical applications it may be convenient to require that as
$\lambda\to\infty$ the numerical solution $f^{n+1}$ and each level
$F^{(i)}$ of the IF method are projected towards the local
Maxwellian without using explicitly the structure of $P(f,f)$. It
is straightforward to verify that this stronger AP property is
satisfied if we replace condition (\ref{eq:cond2}) by \be 0=c_1<
c_2<\ldots< c_\nu < 1. \label{eq:cond3}\ee We remark that the
usual concept of stiff order \cite{MZ} for exponential methods is
in contrast with the latter strong AP property. For example for a
stiff order one method the condition \cite{EX2}
\[
\sum_{i=1}^\nu W_i(\lambda) = \frac{1-e^{-\lambda}}{\lambda}
\]
implies that the Maxwellian term in (\ref{eq:rk2}) vanishes. So
classical stiff order two ETD methods do not satisfy the AP
property (see \cite{MZ} for example).

We conclude this section with a results concerning an important
convexity property of the schemes.

\begin{theorem}
If an explicit exponential Runge-Kutta method in the form
(\ref{eq:rk1})-(\ref{eq:rk2}) satisfies
\begin{eqnarray}
\sum_{j=1}^{i-1}A_{ij}(\lambda)&\leq&
\frac{1-e^{-c_i\lambda}}{\lambda},\quad \forall\,\,\lambda \geq 0,
\quad
i=1,\ldots,\nu \label{eq:cond4}\\
\sum_{i=1}^{\nu}W_{i}(\lambda)&\leq&
\frac{1-e^{-\lambda}}{\lambda},\quad \forall\,\,\lambda \geq 0,
\label{eq:cond5}\end{eqnarray} with $A_{ij}(\lambda)\geq 0$ and
$W_i(\lambda)\geq 0$ then it is unconditionally positive and
entropic.
\end{theorem}

In fact, the operator $P(f,f)$ is nonnegative and convexity of the
schemes is guaranteed if (\ref{eq:cond4}) and (\ref{eq:cond5})
hold. Moreover since $H(M)\leq H(f^n)$ and (see \cite{CV}) \be
H\left(\frac{P(f,f)}{\mu}\right)\leq H(f), \ee by convexity we
have also
\begin{eqnarray*}
\nonumber H(F^{(i)})&\leq &e^{-c_i\lambda}H(f^n)+\lambda
\sum_{j=1}^{i-1}A_{ij}(\lambda)H(F^{(j)})\\
\\[-.25cm]
\nonumber &+& \left(1-e^{-c_i\lambda}-\lambda
\sum_{j=1}^{i-1}A_{ij}(\lambda)\right)H(M)\\&\leq& H(f^n), \qquad
i=1,\ldots,\nu\\
\nonumber H(f^{n+1})&\leq
&e^{-\lambda}H(f^n)+\lambda\sum_{i=1}^{\nu}W_i(\lambda)H(F^{(i)})\\
\\[-.25cm]
\nonumber
&+&\left(1-e^{-\lambda}-\lambda\sum_{i=1}^{\nu}W_i(\lambda)\right)H(M)\\&\leq
&H(f^n).
\end{eqnarray*}

For convexity of an IF method we require
\begin{eqnarray*}
\sum_{j=1}^{i-1}a_{ij}e^{c_j\lambda}&\leq&
\frac{e^{c_i\lambda}-1}{\lambda},\quad \forall\,\,\lambda \geq 0,
\quad
i=1,\ldots,\nu\\
\sum_{i=1}^{\nu}w_{i}e^{c_i\lambda}&\leq&
\frac{e^{\lambda}-1}{\lambda},\quad \forall\,\,\lambda \geq 0.
\end{eqnarray*}
By Taylor expansion we obtain conditions
\begin{eqnarray}
\sum_{j=1}^{i-1}a_{ij}{c_j}^k &\leq& \frac{c_i^k}{k+1},\quad
k=0,1,2,\ldots, \quad i=1,\ldots,\nu \label{eq:cond6}
\\
\sum_{i=1}^{\nu}w_{i}c_i^k&\leq& \frac{1}{k+1},\quad
k=0,1,2,\ldots,\label{eq:cond7}
\end{eqnarray}
This allows to state the following.
\begin{proposition}
An explicit IF method is unconditionally positive and entropic if
the underlying Runge-Kutta method has nonnegative coefficients and
weights satisfying (\ref{eq:cond6})-(\ref{eq:cond7}). \label{pr:6}
\end{proposition}

Note that the above conditions on the choice of the underlying
method are quite restrictive and that we are not using the
bilinearity of $P(f,f)$ which would lead to weaker constraints on
$a_{ij}$ and $w_i$. For example, if we consider the family of
second order Runge-Kutta methods with two levels characterized by
$w_1=1-w$, $w_2=w$ and $a_{21}=(2w)^{-1}$ with $w\in [0,1]$,
Proposition \ref{pr:6} is satisfied when
\[
w^{k-1} \geq \frac{k+1}{2^k},\quad k=0,1,\ldots
\]
which implies $w\in[\frac34,1]$. Examples of methods that satisfy
convexity are the second order Midpoint or Runge method ($w=1$)
and the third order Heun method but not the classical fourth order
Runge-Kutta method \cite{HairerWanner}.

\begin{remark}
Convexity is an essential property if one wants to rely on Monte
Carlo techniques for the computation of the collision operator. In
fact, the resulting scheme is a convex combination of probability
densities and then can be evaluated using the same methods
described in \cite{PR}.
\end{remark}

\subsection{Generalizations and implementation}
\label{ss:im} An essential aspect in the reformulation of the
problem given by (\ref{eq:RIF}) is the choice of the value of the
constant $\mu$ used in estimating the spectrum of the collision
operator. Of course such constant can be chosen at each time step in
order to improve our estimate. In the sequel we show different
choices in the case of variable hard spheres.

The choice of an upper bound for the loss part of the collision term
leads to take $\mu=\mu_p$ where \be
\mu_{p}=\sup_{v}\int_{\R^{3}}|v-v_*|^{\gamma} f(v_*)\,dv_*.
\label{eq:mup}\ee Positivity is guaranteed since this choice implies
clearly $P(f,f)\geq 0$. From a practical viewpoint computation of
(\ref{eq:mup}) can be done at $O(N\log N)$ for a deterministic
method based on $N$ parameters for representing $f(v)$ on a regular
mesh. This can be done using the FFT algorithm thanks to the
convolution structure of the loss term in (\ref{eq:mup}). Thus, in
general, the computation of $\mu_{p}$ will not affect the
computational cost of the scheme.

For Monte Carlo methods based on $v_1,\ldots,v_N$ particles one
should estimate \[ \mu_{p} \approx
\max_{v_i}\frac1{N}\sum_{j}|v_i-v_j|^{\gamma}. \] To avoid the
$O(N^2)$ cost it is a common choice to consider the following upper
bound \be {\tilde \mu}_{p}=2^{\gamma}\max_{i} |v_i-u|^{\gamma}\geq
\max_{v_i}\frac1{N}\sum_{j}|v_i-v_j|^{\gamma}, \quad
u=\frac1{N}\sum_{i} v_i. \label{eq:mc1}\ee

However, such positivity constraint on $P(f,f)$ typically leads to
overestimates of the true spectrum of the collision operator,
especially in Monte Carlo simulations. A better estimate of $\mu$
would be given by the average collision frequency \be
\mu_{a}=\int_{\R^3}\int_{\R^{3}}|v-v_*|^{\gamma}
f(v)f(v_*)\,dv_*\,dv. \label{eq:mua}\ee This can be computed again
at $O(N\log N)$ cost in a deterministic setting whereas in a Monte
Carlo simulation we have
\[
\mu_a\approx \frac{1}{N^2} \sum_{i,j} |v_i-v_j|^{\gamma},
\]
which, to avoid the quadratic cost, can be overestimated by \be
{\tilde \mu}_{a}=\frac{2^{\gamma}}{N}\sum_{i} |v_i-u|^{\gamma}.
\ee

Finally, as suggested in \cite{Filbet}, $\mu$ can be chosen as an
estimate of the spectral radius of the linearized operator $Q$
around the Maxwellian $M$. In fact
\[
Q(f,f)\approx Q(M,M)+\nabla Q(M,M)(M-f)=\nabla Q(M)(M-f),
\]
where $\nabla Q(M,M)$ is the Frechet derivative of $Q$ evaluated
at $M$. For example one can take \be \mu_s = \sup_v
\left|\frac{Q(f,f)}{f-M}\right|. \label{eq:mus}\ee Note however
that the above estimate can be computed easily only in a
deterministic framework.

The choices $\mu=\mu_a$ or $\mu=\mu_s$, although more accurate, pose
the question of stability of the resulting scheme since they do not
guarantee $P(f,f) \geq 0$. Note that, if we denote by
$P_p(f,f)=Q(f,f)+\mu_p f$ and by$\lambda_p=\mu_p\Delta
t/\varepsilon$, for an arbitrary $\mu\geq 0$, we can rewrite
(\ref{eq:rk1})-(\ref{eq:rk2}) as
\begin{eqnarray}
\nonumber
F^{(i)}&=&e^{-c_i\lambda}f^n-(\lambda_p-\lambda)\sum_{j=1}^{i-1}A_{ij}(\lambda)
F^{(j)}+ \lambda_p \sum_{j=1}^{i-1}A_{ij}(\lambda)\frac{P_p(F^{(j)},F^{(j)})}{\mu_p}\\
\label{eq:rk1p}
\\[-.25cm]
\nonumber &+&
\left(1-e^{-c_i\lambda}-\lambda\sum_{j=1}^{i-1}A_{ij}(\lambda)\right)M,
\qquad
i=1,\ldots,\nu\\
\nonumber
f^{n+1}&=&e^{-\lambda}f^n-(\lambda_p-\lambda)\sum_{i=1}^{\nu}W_{i}(\lambda)
F^{(j)}+\lambda_p\sum_{i=1}^{\nu}W_i(\lambda)\frac{P_p(F^{(i)},F^{(i)})}{\mu_p}\\
\label{eq:rk2p}
\\[-.25cm]
\nonumber
&+&\left(1-e^{-\lambda}-\lambda\sum_{i=1}^{\nu}W_i(\lambda)\right)M.
\end{eqnarray}
For stability now we must perform the same analysis of Section
\ref{sec:3-3} on system (\ref{eq:rk1p})-(\ref{eq:rk2p}). For
brevity here we omit the resulting conditions which typically
cannot be satisfied without introducing a stability restriction on
the time step.

To illustrate this let us consider the case $\nu=1$ for IF
methods. We have  \be
f^{n+1}=(1-\lambda_p+\lambda) e^{-\lambda} f^{n}+\lambda_p
e^{-\lambda}\frac{P_p(f^n,f^n)}{\mu_p}+\left(1-e^{-\lambda}-\lambda
e^{-\lambda}\right)M. \label{eq:if1bm} \ee Convexity is guaranteed
as soon as $\Delta t$ satisfies
\[
\Delta t \leq \frac{\varepsilon}{\mu_p-\mu},
\]
or equivalently if we take $\mu\geq\mu_p-\varepsilon/\Delta t$
which represents the lower bound for $\mu$ that makes the scheme
unconditionally positive. On the other hand contractivity and
asymptotic stability remain valid under the weaker restriction
\[
\frac{\mu}2+\frac{\varepsilon}{\Delta t}(1+e^{\mu\Delta
t/\varepsilon})\geq \mu_p.
\]
Similar considerations hold for higher order IF schemes. In such
cases the $AP$ property is guaranteed provided that the underlying
method satisfy Proposition \ref{Pr:4}.

\begin{remark}
Along this paragraph we have assumed $\mu$ constant during the
time stepping. In practice it is clear that the computation of
$\mu_p$ from (\ref{eq:mup}) on the initial data does not guarantee
positivity of all terms in the vector ${\bar f}_1$. In principle
one can take $\mu=\mu(t)$ and rewrite the exponential methods for
a time dependent $\mu$ from \be \frac{\partial
(f-M)e^{\frac1{\varepsilon}\int_0^{t}\mu(s)\,ds}}{\partial
t}=\frac1{\varepsilon}(P(f,f)-\mu(t)
M)e^{\frac1{\varepsilon}\int_0^{t}\mu(s)\,ds}, \label{eq:15b} \ee
and then recompute $\mu_p$ at each stage of the Runge-Kutta level.
This procedure however may be quite expensive and in practice
violation of positivity is rarely observed. In such circumstances
one can set initially a numerical tolerance on $\mu_p$ or simply
repeat the computation with a larger value of $\mu_p$.
\end{remark}

\section{Numerical results}
In this section we perform some numerical tests for the
exponential Runge Kutta schemes applied to the case of the full
Boltzmann equation. We consider the first order IF scheme
(\ref{eq:if1b}) and the second and third order IF schemes obtained
using second order Midpoint and third order Heun methods
respectively. All methods satisfy Proposition \ref{pr:6} and thus
can be implemented using Monte Carlo strategies for the evaluation
of the five fold collision integral \cite{PR}. Since we are
interested in measuring the time discretization error of the
different schemes we need to cancel the other sources of error in
the computations. This is realized using a very large number of
particles and statistical averages.

We consider two different test cases: first the evolution of the
fourth order moment in a space homogeneous case
and then the heat flux behavior in a space inhomogeneous shock
problem.

\subsection{Homogeneous relaxation}

As initial data we consider an equilibrium Maxwellian distribution
with temperature $T=6$, density $\varrho=1$ and mean velocity
$u=-0.5$. To this distribution we add a bump on the right tail along
the x-axis. The bump is realized adding a Maxwellian with mass
$\varrho_{b}=0.5\ \varrho$, mean velocity $u_{b}=4 \ \sqrt{T}$ and
temperature $T_{b}=0.5 \ T$ to the initial Maxwellian distribution.
We consider the case of Maxwell interactions and hard spheres. The
simulations are run till the equilibrium is approximately reached,
which means $t=0.4$ in the case of hard spheres and $t=0.8$ for
Maxwell molecules. The reference solution is computed by the Bird
method \cite{bird} which converges toward the exact solution when
the number of particles goes to infinity.

\begin{figure}
\begin{center}
\includegraphics[scale=0.4]{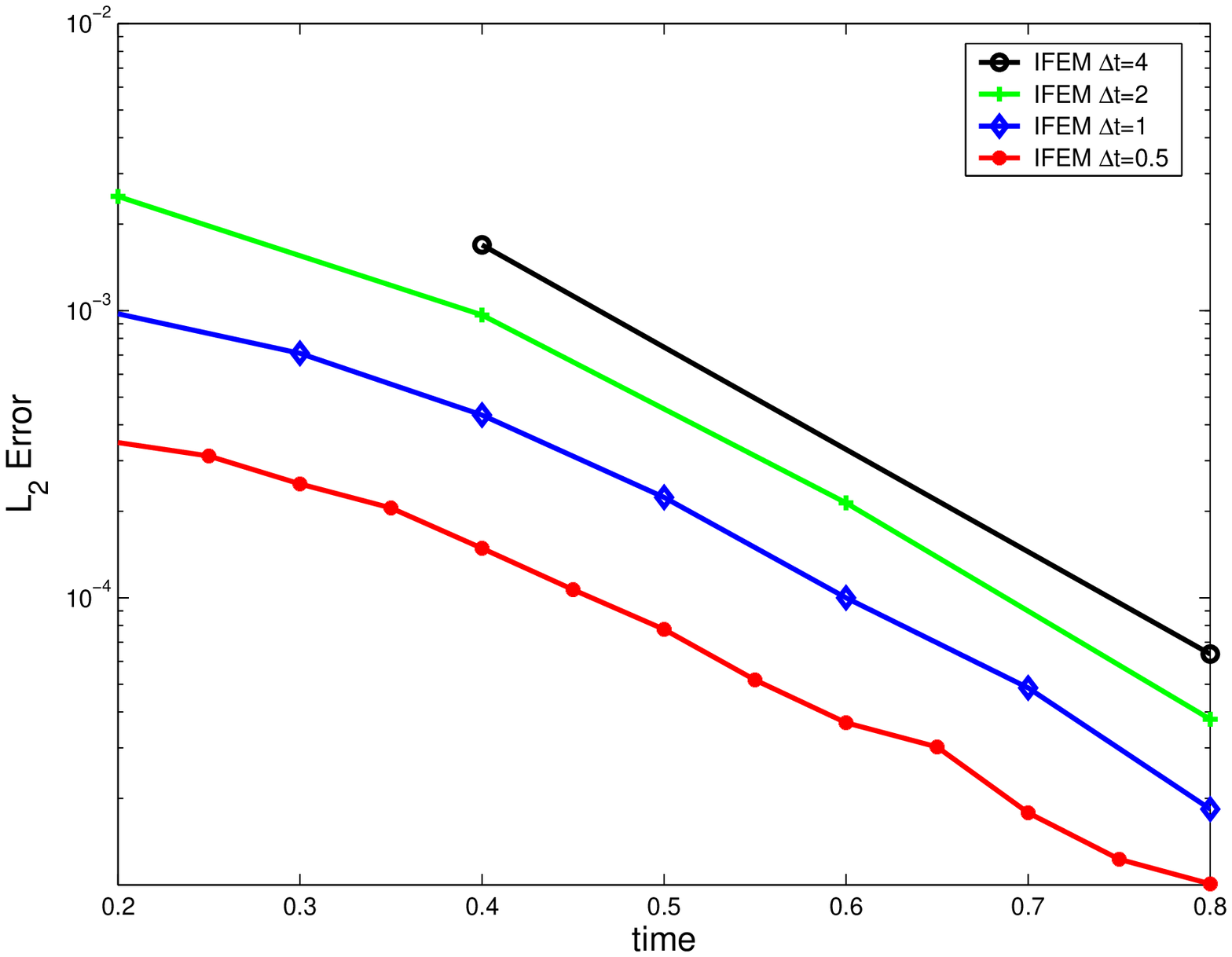}
\includegraphics[scale=0.4]{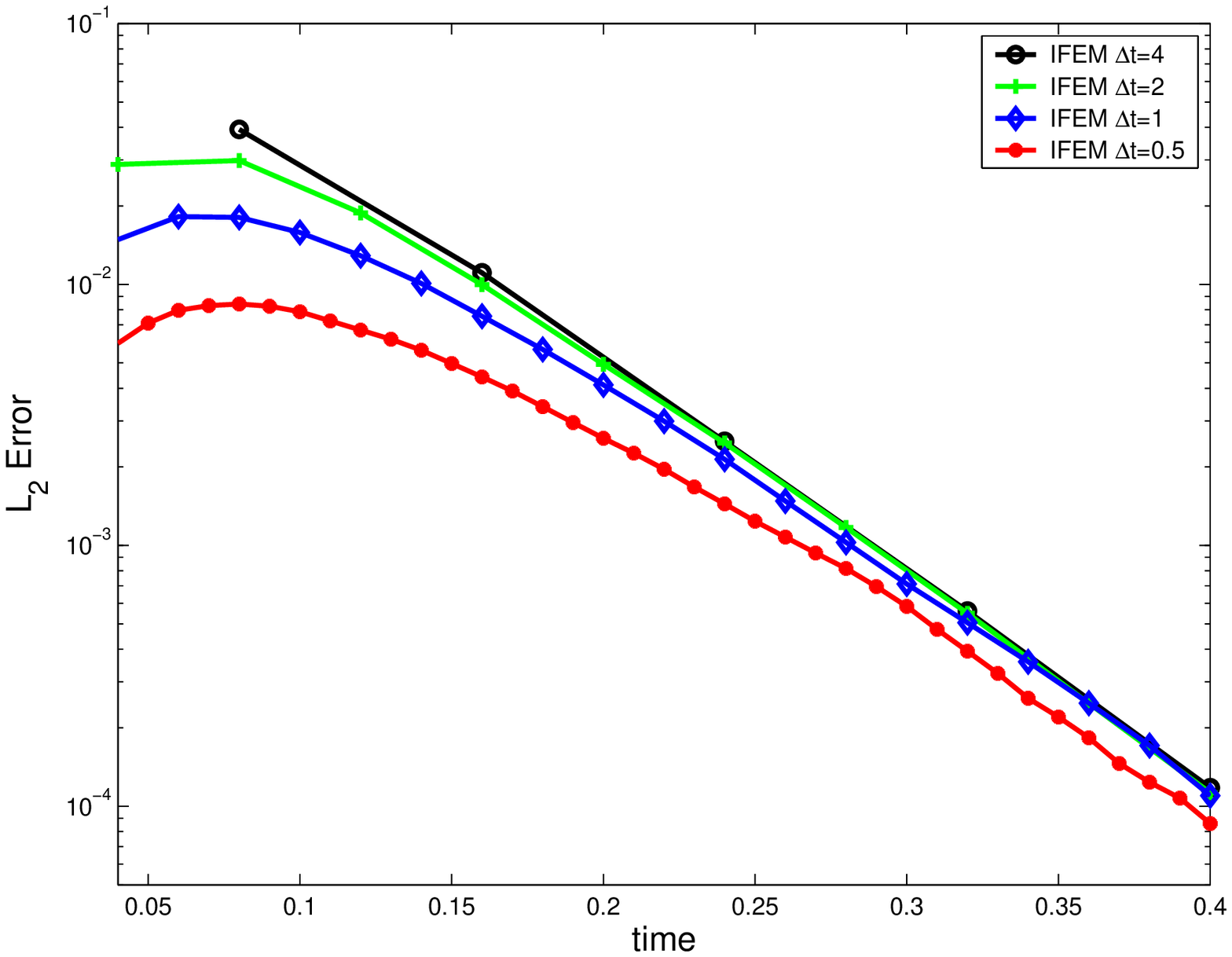}
\includegraphics[scale=0.4]{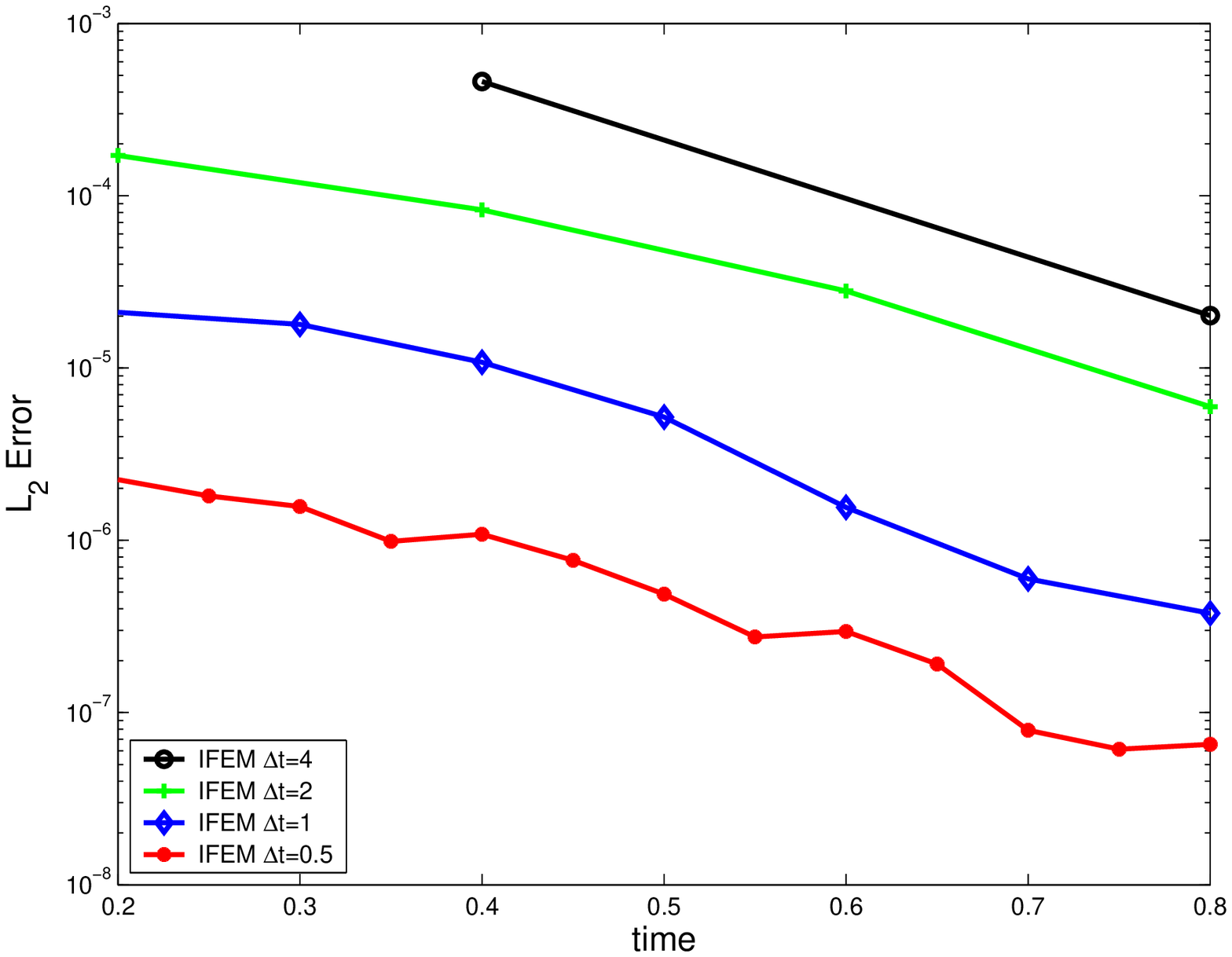}
\includegraphics[scale=0.4]{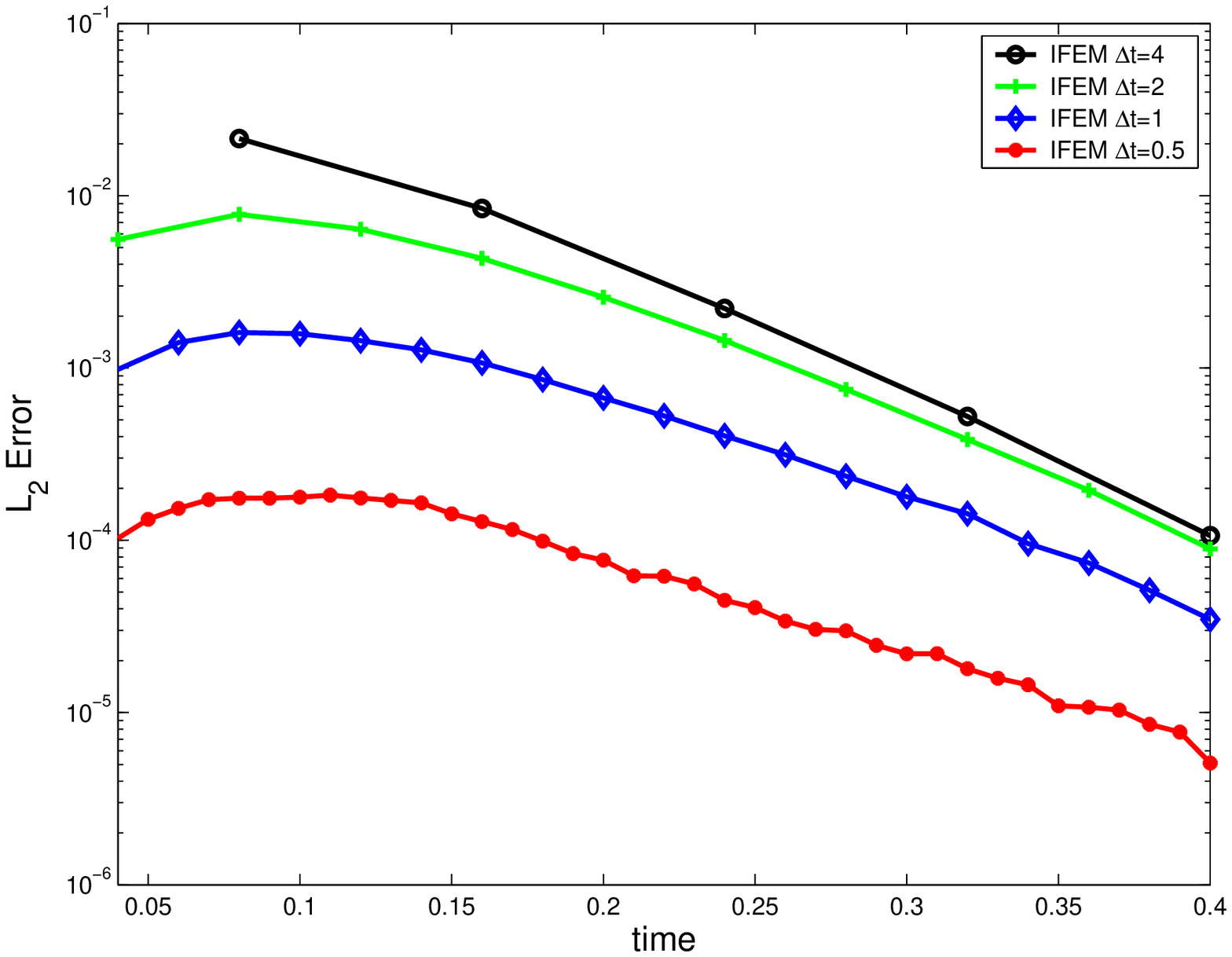}
\includegraphics[scale=0.4]{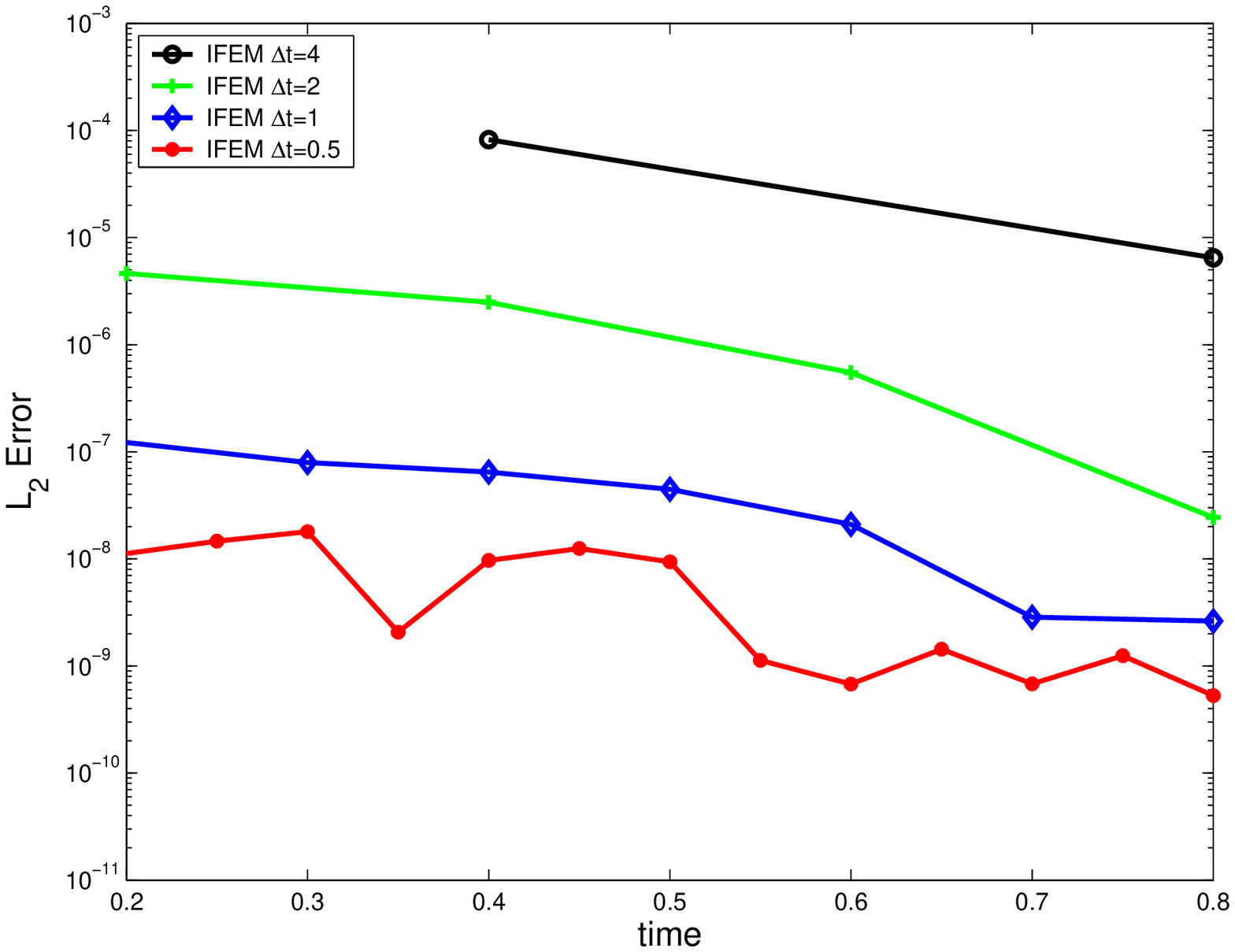}
\includegraphics[scale=0.4]{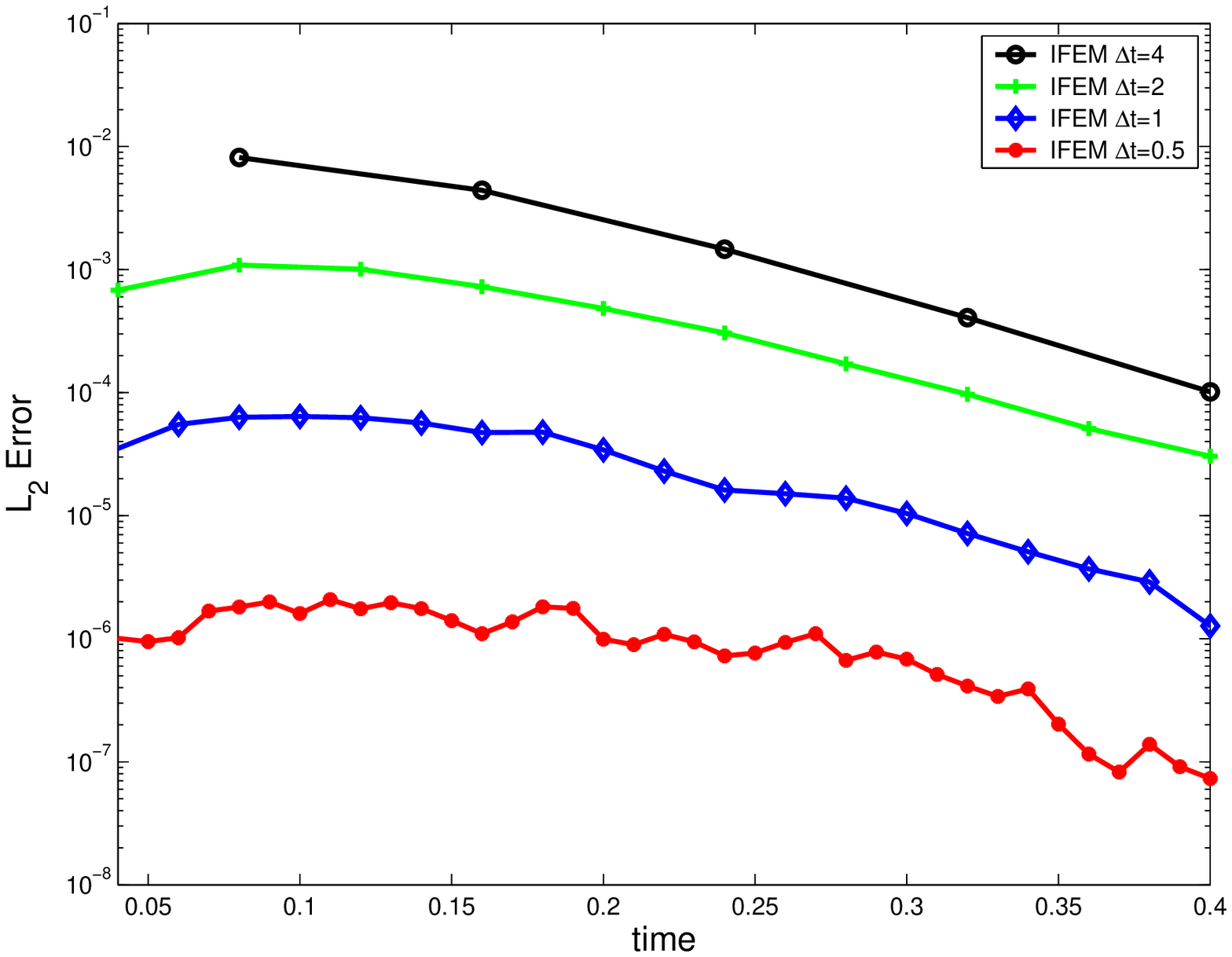}
\caption{$L_{2}$ Error for the Fourth Moment Relaxation for the
homogeneous relaxation problem with Maxwellian
particles (left) and hard spheres (right).}\label{fig:DSMC1}
\end{center}
\end{figure}

In Figure \ref{fig:DSMC1} we show the $L_{2}$ error for the fourth
order moment of the distribution function $f$ for Maxwellian
molecules and hard spheres. In each of the plots the error is
depicted for different choices of the time step: respectively $\mu
\Delta t/\varepsilon=0.5, 1, 2$ and $4$. In the case of Maxwell
molecules $\mu=1$ while for hard spheres $\mu\gg 1$ is an upper
bound for the collision cross section computed using the simple
choice (\ref{eq:mc1}). As a consequence for the same values of
$\mu \Delta t /\varepsilon$ the $L_2$ norm of the error is larger
for hard spheres with respect to Maxwell molecules. Here we do not
use any strategy to reduce the effect of possible overestimation
of $\mu$ as described in section \ref{ss:im}. We leave this
possibility to further research. The expected convergence rate is
observed for all schemes.


\subsection{A shock wave computation}
We consider a Sod shock tube test with initial values
\begin{equation}
\nonumber
\begin{array}{l}
\textbf{u}_{L}=\left(
\begin{array}{ll}
1 \\0\\5
\end{array}
\right), \ \hbox{if} \ 0  \leq x <0.5 \ \ \
\textbf{u}_{R}=\left(\begin{array}{l} 0.125 \\ 0 \\ 4
\end{array}\right), \ \hbox{if} \ 0.5\leq x \leq 1.
\end{array}
\label{eq:Sod}
\end{equation}
The solution is computed using $150$ grid points in space, the
final time is $t=0.05$. The transport step is solved exactly by
particle transport as it is usual in Monte Carlo methods. As
before we consider a very large number of particles and averaged
the solution over several runs. Again the reference solution is
obtained letting the time step go to zero and the number of
particles to infinity using Bird's DSMC method.

\begin{figure}
\begin{center}
\includegraphics[scale=0.4]{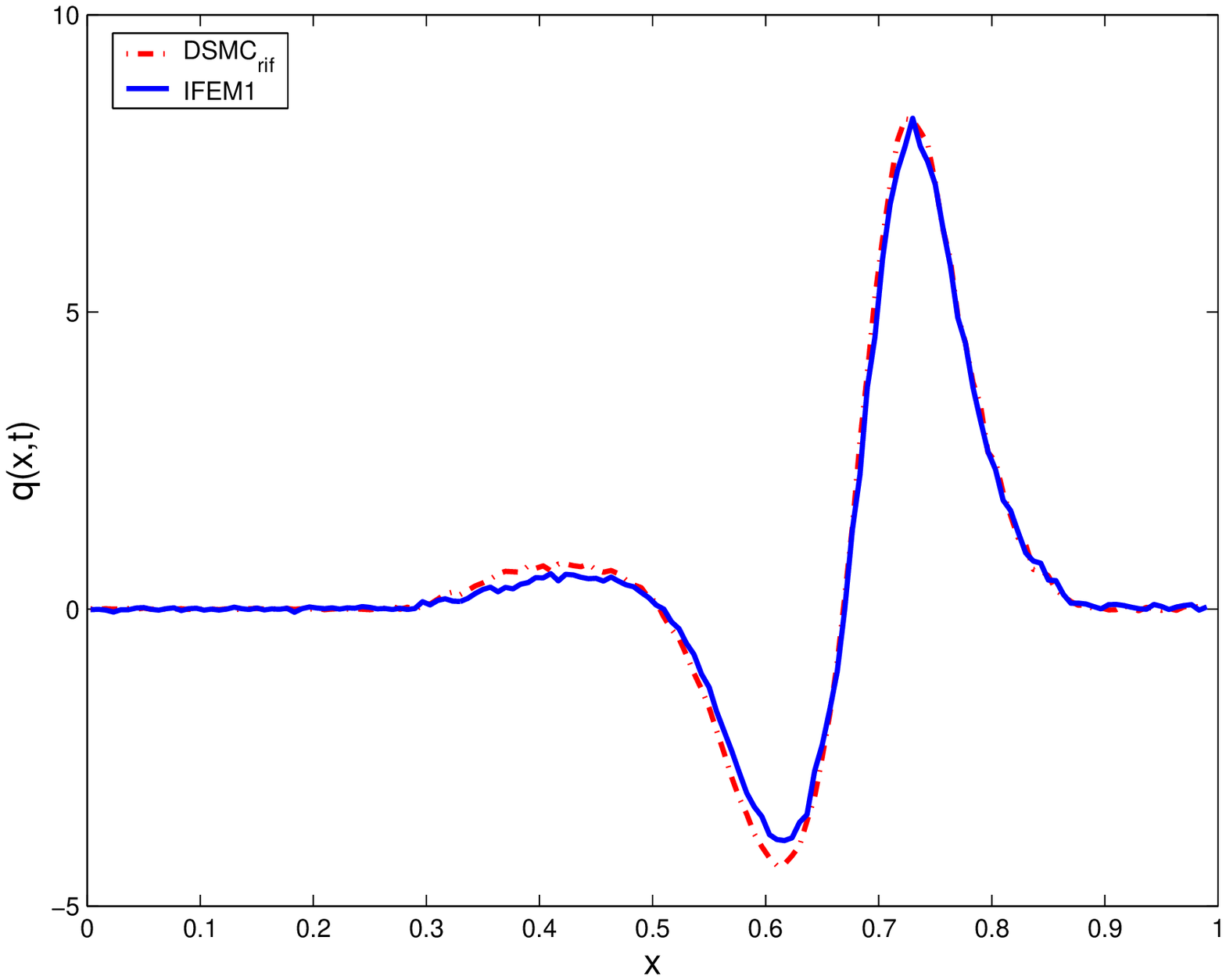}
\includegraphics[scale=0.4]{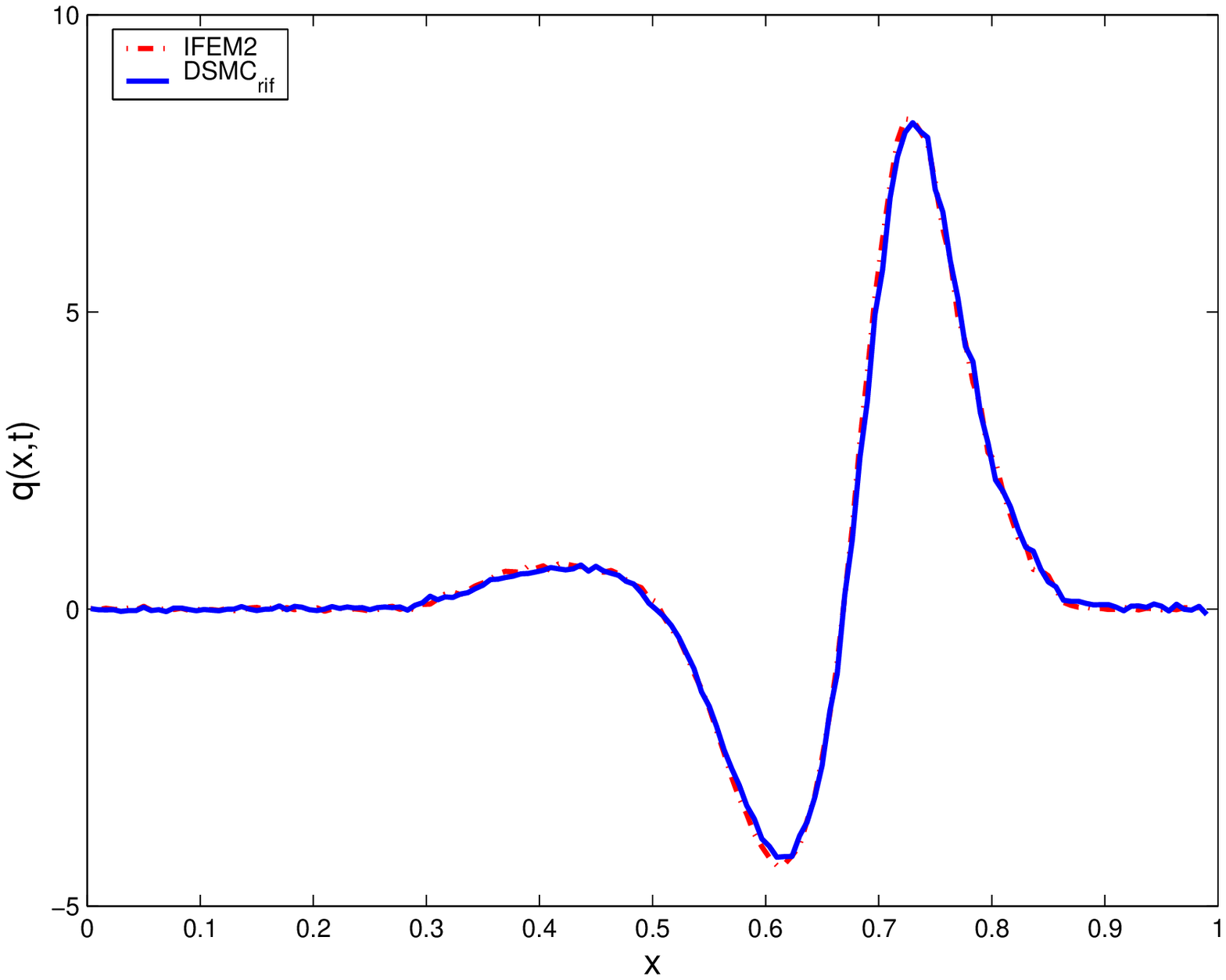}\\
\includegraphics[scale=0.4]{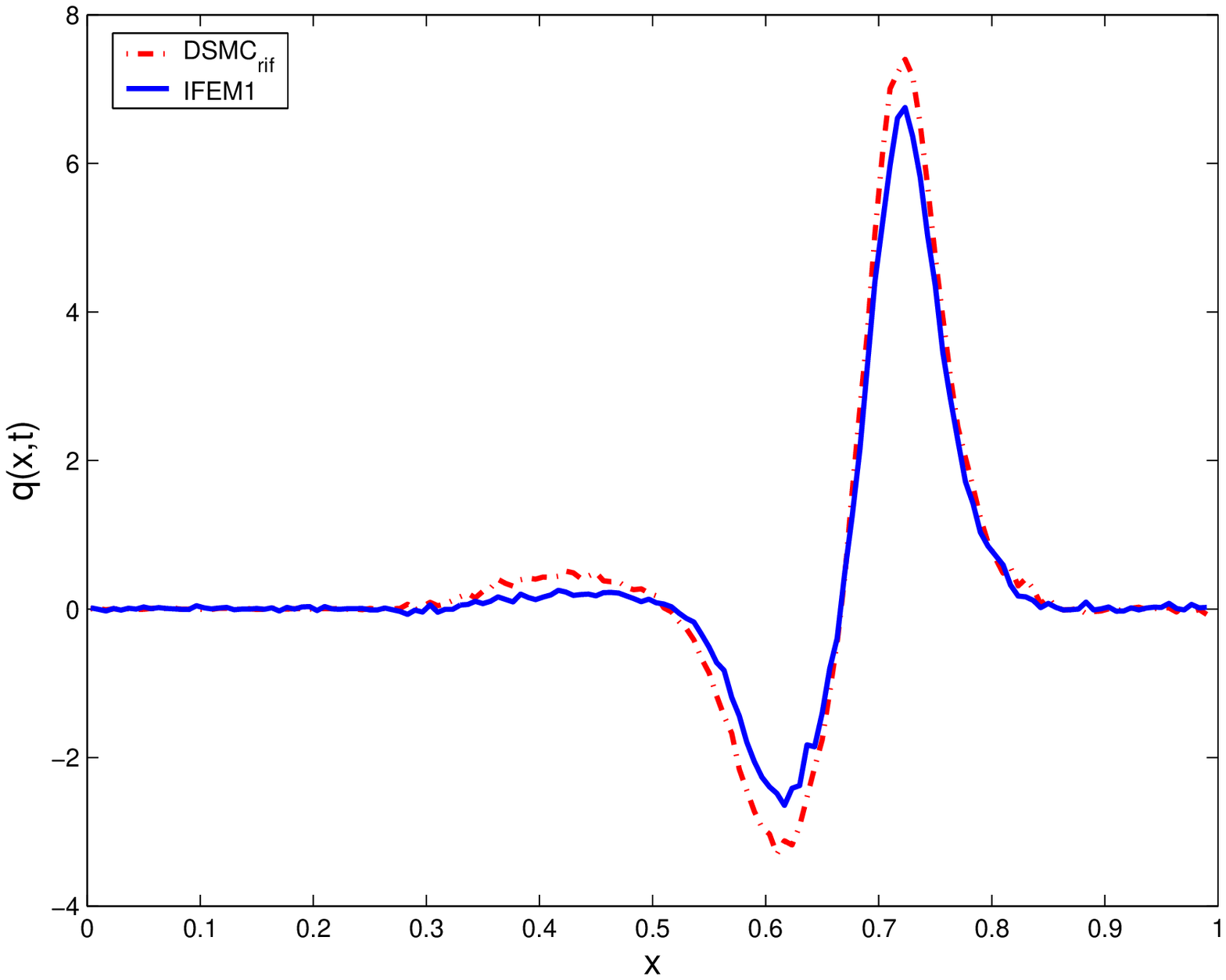}
\includegraphics[scale=0.4]{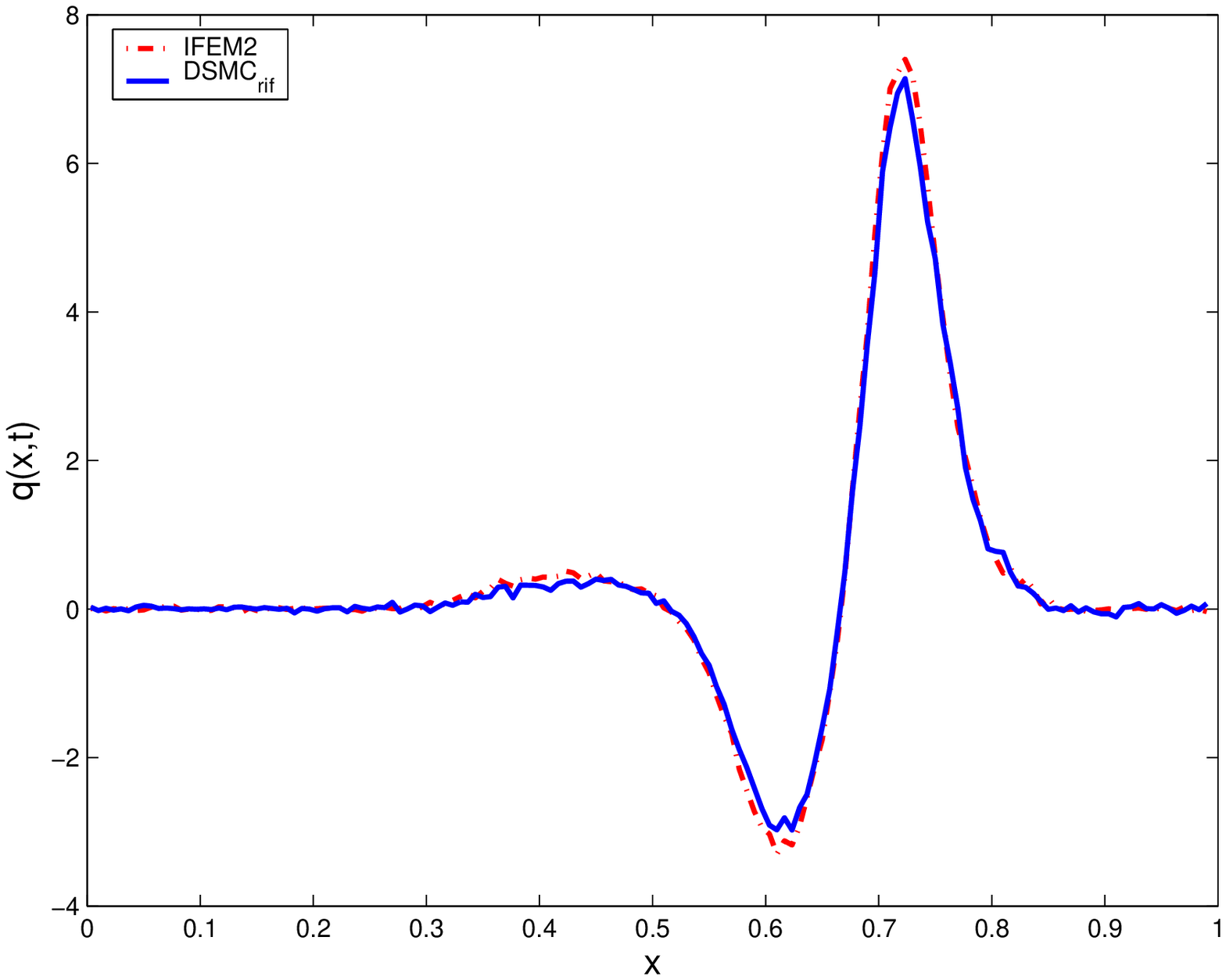}\\
\includegraphics[scale=0.4]{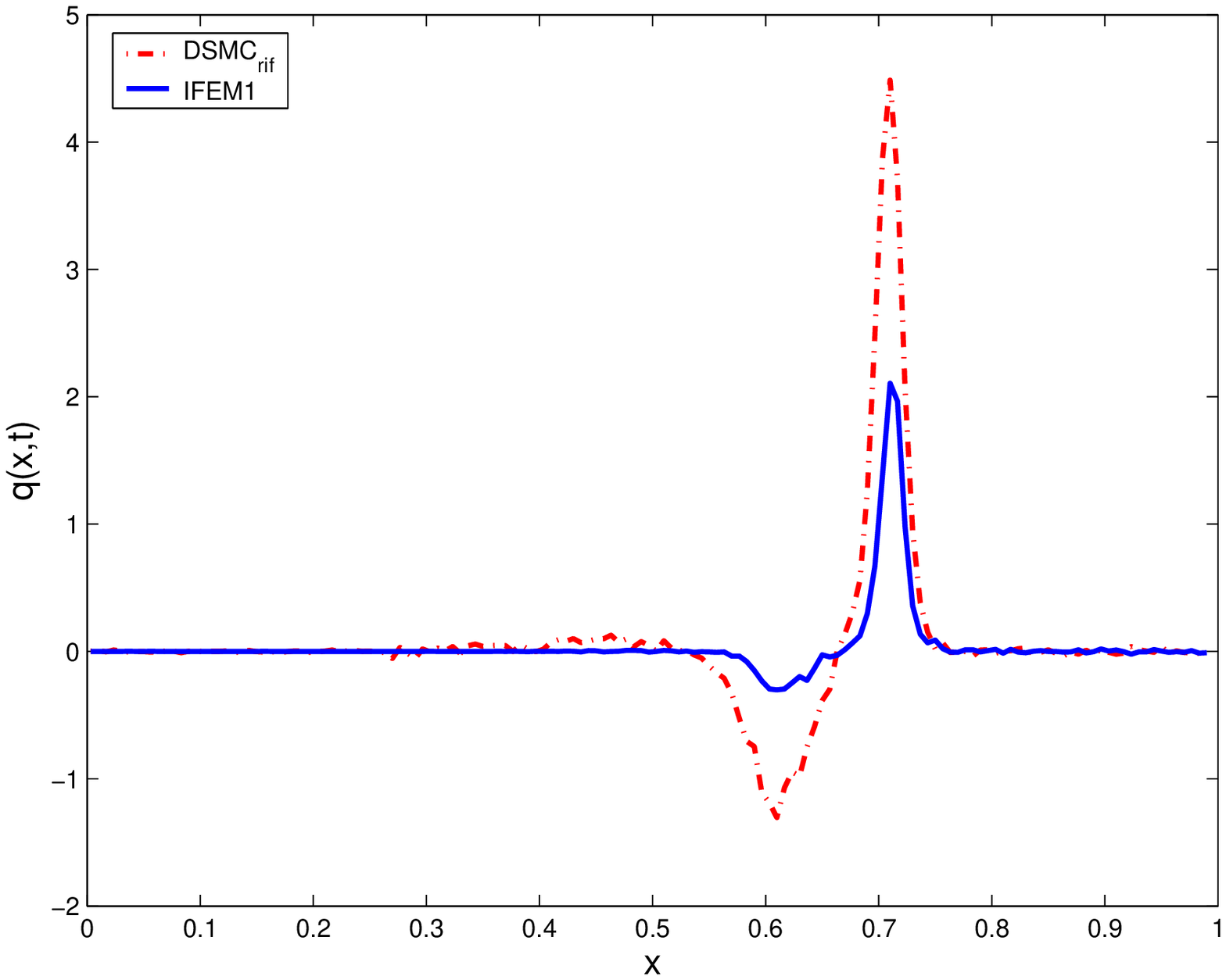}
\includegraphics[scale=0.4]{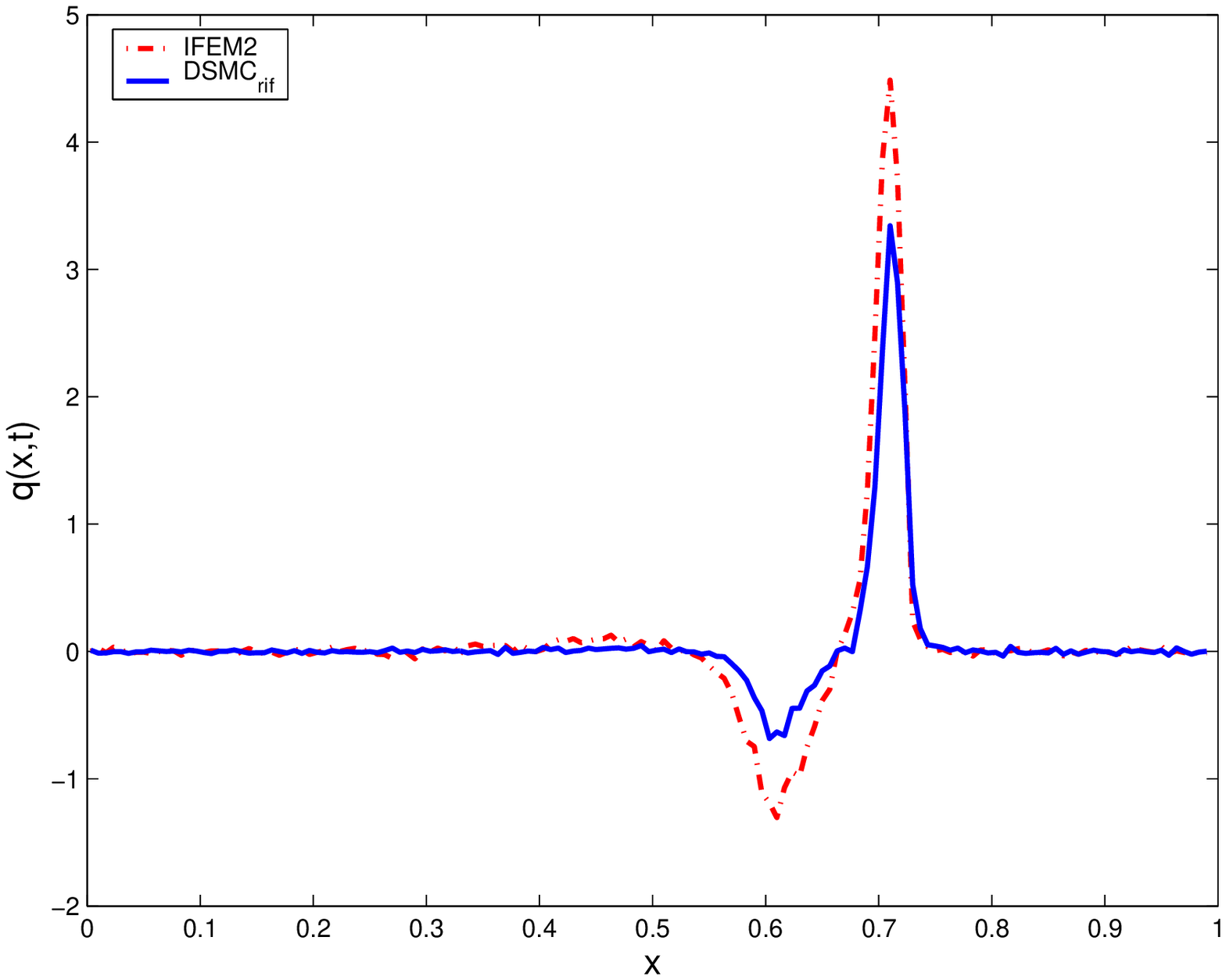}
\caption{Heat flux profile for first order (left) and second order (right) IF schemes. Top
Kundsen number $\varepsilon=10^{-3}$, middle $\varepsilon=5 \
10^{-4}$ and bottom $\varepsilon=10^{-4}$.
$\Delta t=10^{-3}$.}\label{fig:SOD1}
\end{center}
\end{figure}

We report the heat flux profile for the first and the second order
Runge Kutta integration factor scheme in figures \ref{fig:SOD1}. A
first order splitting is employed for the first order IF while a
second order Strang splitting is used in combination with the
second order IF method. From top to bottom of the figures the
Knudsen number values are $\varepsilon=10^{-3}, 5 \ 10^{-4}$ and
$10^{-4}$. The time step is fixed equal to $10^{-3}$.

Both schemes show a good agreement with the reference solution for
$\varepsilon=10^{-3}$. Then, when the Knudsen number is halved we
start to see some discrepancies between the profiles of the first
order method and the reference solution. On the contrary the
second order method is still in good agreement with the reference
solution. When the Knudsen number becomes $1/10$ of the time step
both the schemes present larger errors but they are still able to
catch correctly the departure from equilibrium. This is especially
true for the second order scheme which is able to reproduce the
correct profile with sufficient accuracy.


\section{Conclusions and further developments}

We have presented a general class of exponential schemes for the
numerical solution of nonlinear kinetic equations in stiff regimes.
The schemes generalize the class of schemes previously developed in
\cite{toscani} and share the fundamental property of asymptotic
preservation. Even if the schemes have been developed in the case of
the Boltzmann equation for dilute gases, extension of the schemes to
other collisional kinetic equations which possess a smooth
equilibrium state are straightforward. Let us also mention that
decomposition (\ref{eq:dec}) represents only one of the possible
choices in order to linearize the collision operator close to
equilibrium and then using it as a penalization factor in the
construction of the numerical methods. For example a more accurate
penalization can be obtained using the so called {\it ES-BGK}
equilibrium function \cite{Hol}, instead of the standard Maxwellian,
which is well known to provide a better approximation of the
collision dynamics. Let us finally mention that in principle the
schemes can be extended to the non-splitting method case for
(\ref{eq:1}) using the integral representation
\[
f(x,v,t)=f(x,v,0)e^{-\mu
t/\varepsilon}+\frac{\mu}{\varepsilon}\int_{0}^t e^{-\mu
(t-s)/\varepsilon} G(y(s))\,ds+\frac{\mu}{\varepsilon}\int_{0}^t
e^{-\mu (t-s)/\varepsilon} M(x,v,s)\,ds,
\]
where
\[
G(f,s)=-\frac{\varepsilon}{\mu}v\nabla_x f+\frac{P(f,f)}{\mu}- M.
\]
Here we do not explore further this direction and we leave it to
future researches.




\end{document}